\documentclass[a4paper,10pt]{article}
\usepackage{amsfonts,amsmath,amssymb,enumitem}

\newcommand{\BB}{\mathcal{B}}
\newcommand{\CC}{\mathcal{C}}
\newcommand{\DD}{\mathcal{D}}
\newcommand{\EE}{\mathcal{E}}

\newcommand{\W}{\mathbb{W}}

\newcommand{\HH}{\mathcal{H}}
\newcommand{\OmT}{\Omega^{\eps}_{0,T}}
\newcommand{\R}{\mathbb{R}}
\newcommand{\N}{\mathbb{N}}

\newcommand{\eps}{\varepsilon}
\newcommand{\om}{\omega}

\newcommand{\bsob}{\mathrm{H}_0}
\newcommand{\sob}{\mathrm{H}}
\newcommand{\Sa}{\mathrm{Sa}}

\newtheorem{Lemma}{Lemma}
\newtheorem{Proposition}{Proposition}
\newtheorem{Corollary}{Corollary}

\newtheorem{Theorem}{Theorem}

\begin{document}

\title{A Sub-Gaussian estimate for Dirichlet Heat Kernels on Tubular Neighbourhoods and Tightness of Conditional Brownian Motion}

\author{Olaf Wittich\\Lehrstuhl I f\"ur Mathematik\\RWTH Aachen University}

\maketitle

\abstract{We prove tightness of a family of path measures $\nu_{\eps}$ on tubes $L(\eps)$ of small diameters around a closed and connected submanifold $L$ of another Riemannian manifold $M$. Together with a result from \cite{NobWit:19}, that implies weak convergence of the measures as the tube radius $\eps$ tends to zero to a measure supported by the path space of the submanifold. As a consequence, we obtain weak convergence of the measures obtained by conditioning Brownian motion to stay within the tubes $L(\eps)$ up to a finite time $T>0$, and we identify the limit measure.\\


}
	
\section{Introduction}\label{intro} We consider an $l$-dimensional closed and connected submanifold $L\subset M$ of a Riemannian manifold $M$ of dimension $m > l$. Let $r >1$ be a number such that the $r$-tube $L(r)$ is mapped diffeomorphically onto the $r$-neighbourhood $NL(r)$ of the zero section in the normal bundle by the exponential map. For $0<\eps < r$, we consider Brownian motion conditioned to remain within $L(\eps)$ up to a finite time horizon $T>0$. In \cite{NobWit:19}, we showed that the marginals of this family of conditional processes tend to the marginals of a process supported by the path space of the submanifold with a path measure which is equivalent to the Wiener measure on $L$. In this paper, we prove that the family of conditional processes is actually tight, showing that the path measures of the conditional processes actually converge weakly to the limit measure. For embeddings into Euclidean space this was shown in \cite{SidSmoWeiWit:03b}. \\

\noindent As in \cite{NobWit:19}, we do not attempt to give a direct proof for conditional Brownian motion in the first place, but consider measures $\nu_{\eps}$, equivalent to the measure of Brownian motion absorbed at the tube boundary $\partial L(\eps)$. These measures are related to the generators $H_{\eps}$ constructed in \cite{NobWit:19} in such a way, that the associated semigroups correspond to the one-dimensional marginals of $\nu_{\eps}$. By the same rescaling and renormalization procedure as in \cite{NobWit:19}, we now obtain ultracontractive symmetric Markov semigroups on $L^2(L(1), \mu_{\Sa})$ given by smooth kernels for $t>0$. Here $\mu_{\Sa}$ denotes the Riemannian volume associated to the Sasaki-metric on $L(1)$. Using properties of the ground states of the generators from Proposition \ref{asymptotics} and Corollary \ref{unifestground}, together with a Hardy inequality, we prove in Proposition \ref{Rosen} a version of the Rosen Lemma by which the logarithmic Sobolev inequalities of the generator can be used to obtain intrinsic ultracontractivity of the semigroups. The corresponding $L^2$-$L^{\infty}$ operator bound  can be extended to the sub-Gaussian estimate Proposition \ref{preparation} of the kernel, by methods described for instance in \cite{Dav:89}, Section 3.2. From the uniform estimate of the kernel, we conclude tightness of the measure family $(\nu_{\eps}:\eps>0)$ and finally of conditional Brownian motion by the standard argument Proposition \ref{Kolmogorov}. \\ 

\noindent In the first section, we reintroduce the setup from \cite{NobWit:19}, in the second section, we prove by a rather general argument that tightness of the conditional processes is equivalent to tightness of its projections, which are processes on the path space of $L$. In Section \ref{3}, we establish the sub-Gaussian bound Corollary \ref{subgauss} and conclude tightness of the projected processes. The main work, namely the proof of Proposition \ref{preparation}, is done in Section \ref{beweis}. Everything that is done here, follows closely the considerations in \cite{Dav:89}. However, we attempted to outline the basic ideas from this source for the convenience of the reader.

\subsection{The conditional process}\label{setup} In this section, we are going to describe the setup from \cite{NobWit:19}.

\noindent a. Let $\Omega := C(\lbrack 0,\infty), M)$ be the path space and 
$$
\Omega^{\eps}_{s,t}:= \lbrace \om\in\Omega\,:\,\om(u)\in L(\eps), s\leq u \leq t\rbrace.
$$  
Denoting by $\mathbb{W}$ the Wiener measure on $M$, we fix some finite $T>0$ and consider the measure 
$$
\nu (d\om) = \exp\left(\frac{1}{2}\int_0^{T}\varphi U(\om(s)) \,ds\right)\,\mathbb{W}(d\om),
$$
on the path space of $M$. Here, $\varphi\in C^{\infty}(M)$ is smooth with $\varphi\vert_{L(1)} = 1$ and $\varphi\vert_{M\setminus L(r)}=0$, and for some $r>1$, $U\in C^{\infty}(L(r))$ is a potential given by
\begin{equation}\label{potential}
U:=\rho^{1/2}\Delta\rho^{-1/2},
\end{equation}
where $\rho = d\mu/d\mu_{\Sa}$ is the Radon-Nikodym density of the Riemannian volumes $\mu$ associated to the induced metric and $\mu_{\Sa}$ associated to the Sasaki metric on the tube. Finally, $\Delta$ denotes the Laplace - Beltrami operator on $M$.\\

\noindent Now, we consider the probability measures $\mathcal{L}_{\eps}$, $\eps > 0$, which are obtained by restricting $\nu$ to the set $\Omega^{\eps}_{0,T}$ followed by normalization to total mass one. To be precise, 
\begin{equation*}
\mathcal{L}_{\eps} (d\om) := \frac{\nu (d\om\cap \Omega^{\eps}_{0,T})}{\nu(\Omega^{\eps}_{0,T})},
\end{equation*}
which yields a distribution supported by the path space $\Omega^{\eps}_{0,T}$. The processes with distribution $\mathcal{L}_{\eps}$ are denoted by $(x_t^{\eps})_{0\leq t\leq T}$. For $0\leq s < t\leq T$, the transition kernel of $x_t^{\eps}$, given by the conditional probability
$Q_{\eps}(s,x;t,dy) := \mathcal{L}_{\eps} (\om (t)\in dy\,\vert\, \om (s) = x)$, can be written as
\begin{equation}\label{conditional_kernel}
 Q_{\eps}(s,x;t,dy) = \frac{\pi_{T-t}^{\eps}(y)}{\pi_{T-s}^{\eps}(x)}P^{\eps}(s,x;t,dy),
\end{equation}
where
\begin{eqnarray*}
P_{\eps}(s,x;t,dy) &=& \nu (\om(t)\in dy, \Omega_{s,t}^{\eps}\,\vert\,\om (s) = x) \\
&=& \nu(\om(t)\in dy, t < \tau_{\eps}(\om)\,\vert\,\om (s) = x)
\end{eqnarray*}
is the transition kernel of the process absorbed at the boundary, i.e. $\tau_{\eps}$ is the first exit time from $L(\eps)$, and $\pi^{\eps}_{u}(w) :=\int_{L(\eps)}P^{\eps}(0,w;u,dz)$. Since $M$ is smooth, $P^{\eps}$ and therefore $Q^{\eps}$ have a density with respect to the Riemannian volume measure of $M$.\\

\noindent b. By the Feynman-Kac formula, integration with respect to the transition kernel can be represented by
\begin{eqnarray*}
& & \int_{L(\eps)} f(y)P^{\eps}(s,x;t,dy) \\ 
&=& \int_{\Omega} f(\omega(t))\exp\left(\frac{1}{2}\int_s^{t} U(\om(s)) ds\right)\mathbb{W}(\om(t)\in dy, t < \tau_{\eps}(\om)\,\vert\,\om (s) = x)\\
&=& \left(e^{-\frac{t-s}{2}(\Delta-U)_{\eps}}f\right)(x),
\end{eqnarray*}
where $(\Delta - U)_{\eps}$ denotes the second order differential operator with Dirichlet boundary conditions associated to the quadratic form 
$$
f\mapsto  \int_{L(\eps)} df\wedge \star df - \star Uf^2 , \,\, f\in\bsob^1 (L(\eps),\mu) .
$$
We study the limit of this family of operators as $\eps$ tends to zero by a rescaling procedure, followed by a renormalization of the lowest eigenvalue. To explain this, we first have to note that there are two different metrics on $L(1)$, the metric induced by the embedding $L(1)\subset M$ and the Sasaki metric. The Sasaki metric can most easily described by using that the exponential map of $M$, restricted to the normal bundle, yields a diffeomorphism $\exp\vert_{NL(1)}: NL(1)\to L(1)$ from the one-neighbourhood $NL(1) :=\lbrace W\in NL\,:\, \Vert W\Vert_{NL}<1\rbrace$ of the zero section in the normal bundle to the tube $L(1)$. By pulling back the metric suitably, we can identify any Riemannian structure on $L(1)$ and $NL(1)$ which we will do in the sequel without further mentioning. Then, for $W\in NL(1)$, $X,Y\in T_W NL(1)$, the Sasaki metric is given by
$$
\langle X,Y\rangle_{\Sa} = \langle \pi_* X,\pi_*Y\rangle_L + \langle K_WX, K_WY\rangle_{NL},
$$  
where $\pi :NL(1)\to L$ is the projection, and $K_W :T_W NL\to N_{\pi(W)}L$ is the connection map of the normal bundle. Let now $\sigma_{\eps}:NL(\eps)\to NL(1)$ be given by $W\mapsto \eps^{-1}W$. By the identification described above, we will denote the corresponding map $\exp\circ \eps^{-1}\circ\exp^{-1}: L(\eps)\to L(1)$ by the same symbol.  Now, we consider the rescaling map $\Sigma_{\eps}:L^2(L(1),\mu_{\Sa})\to L^2(L(\eps),\mu)$ given by
\begin{equation}\label{rescaling_map}
\Sigma_{\eps}f:=\left(\eps^{m-l}\rho\right)^{-1/2} \sigma_{\eps}^*f.
\end{equation}
By partial integration, it turns out that 
\begin{equation}\label{laplacemainops}
\Sigma_{\eps}^{-1}\circ(\Delta-U)_{\eps}\circ\Sigma_{\eps} = H_{\eps},
\end{equation}
where $H_{\eps}$ with domain $\bsob^1\cap\sob^2(L(1),\mu_{\Sa})$ is self-adjoint and non-negative on $L^2(L(1),\mu_{\Sa})$ associated to the quadratic form
\begin{equation}\label{qform}
q_{\eps} (f) = \int_{L(\eps)} d\Sigma_{\eps}f \wedge \star d\Sigma_{\eps}f =  \int_{L(1)} df \wedge \#_{\eps} df 
\end{equation}
with $\#_{\eps} \alpha = \eps^{l-m}\,(\sigma_{\eps}^*)^{-1}\circ (\rho\,\star)\circ\sigma_{\eps}\alpha$ for one - forms $\alpha$ on $L(1)$, and the domain of the quadratic form is given by $\DD = \bsob^1(L(1),\mu_{\Sa})$.\\

\noindent c. Because the parameter $\eps > 0$ is closely related to the tube radius, the perturbation problem for $H_{\eps}$ is not to be expected to yield a sensible limit as $\eps$ tends to zero. However, if $\lambda_0>0$ is the smallest eigenvalue for the Dirichlet problem on the $(m-l)$-dimensional Euclidean unit ball $B\subset\R^{m-l}$, the semigroups generated by $H_{\eps}^0 := H_{\eps} - \eps^{-2}\lambda_0$ converge strongly to a semigroup on a subspace $E_0\subset L^2(L(1),\mu_{\Sa})$. The orthogonal projection onto the subspace $E_0$, denoted by the same symbol, is obtained as follows:\\

\noindent The eigenspace corresponding to the lowest eigenvalue $\lambda_0>0$ of the Dirichlet Laplacian on the flat unit ball  is one-dimensional generated by a normed eigenfunction $\varphi_0(v) = \widetilde{\varphi_0}(\Vert v\Vert)\geq 0$, which is invariant under rotations. By \cite{Her:60}, Prop. 3.3, all fibers $F_x := \pi^{-1}(x)$, $x\in L$, of the projection $\pi : L(1)\to L$ are isometric when equipped with the metric induced by the Sasaki metric. Hence, they are isometric to the flat unit ball $B\subset \R^{m-l}$ and $\phi_0:L(1)\to \R$, given by $\phi_0(W) := \widetilde{\varphi_0}(d(W,x))$ is a well-defined, smooth function on $L(1)$. Now, integration along the fibre yields a mapping $\langle \cdot, \phi_0\rangle : L^2 (L(1),\mu_{\Sa}) \to L^2(L,\mu_L)$ where $f\in  L^2 (L(1),\mu_{\Sa})$ is mapped to the function $f_b:L\to\R$ given by
$$
x\mapsto  \langle f,\phi_0\rangle_{x} = \int_{\pi^{-1}(x)} f(W)\phi_0(W) d\mu_x(W) 
$$
on $L$. Thus, the orthogonal projection $E_0$ is given by
$$
\left(E_0\right) f (W)= \phi_0(W) \,\left(\langle f,\phi_0\rangle_{\pi (W)}\right) \in L^2 (L(1),\mu_{\Sa}) .
$$
In particular, every function $f\in E_0$ can be uniquely written $f = (f_b\circ\pi)\phi_0$. Now, the main result of \cite{NobWit:19} reads as follows:

\begin{Theorem}\label{Main} Let $u(\eps)_{\eps \geq 0}\subset L^2(L(1),\mu_{\Sa})$ be a strongly continuous family of functions and denote by $\Delta_L$ the Laplace-Beltrami operator on $L$. Then, for all $n\geq 1$, we have
	\begin{equation}\label{result}
	\lim_{\eps\to 0} e^{-\frac{t}{2}H_\eps^0}u(\eps) = E_0\,e^{-\frac{t}{2}\Delta_L} \,E_0u(0) 
	\end{equation}
	uniformly on each compact sub-interval $I\subset (0,\infty)$ in the Sobolev space $\sob^{2n}(L(1),\mu_{\Sa})$.  
\end{Theorem}

\noindent{\bf Remark.} (a) Using the explicit form of the projection, we obtain
\begin{equation}\label{meaning}
E_0\,e^{-\frac{t}{2}\Delta_L} \,E_0f = \phi_0\,\left(e^{-\frac{t}{2}\Delta_L}f_b\right)\circ\pi.
\end{equation}
(b) Theorem \ref{Main} will still hold if $u(\eps)$ is only strongly continuous at $\eps=0$.\\

\noindent d. In \cite{NobWit:19}, we also concluded from Theorem \ref{Main} the following statement that the one-dimensional marginals  of the processes $x_t^{\eps}$ converge to those of Brownian motion $x^0_t$ on $L$. Since they are Markov processes, this statement implies convergence in finite dimensional distributions.

\begin{Corollary}\label{findim_conv} Let $x^0_{t, 0\leq t\leq T}$ be Brownian motion on $L$. Let $x_0^{\eps} = x\in L$ be a fixed common starting point. Then, for all $f\in C^{\infty}(M)$ and $0\leq t\leq T$, we have 
$$
\lim_{\eps\to 0} 
E^x\lbrack f(x_t^{\eps})\rbrack = E^x\lbrack f\vert_L(x_t^0)\rbrack,
$$
i.e. the associated flows converge as $\eps$ tends to zero. 
\end{Corollary}

\noindent Convergence of all marginals is the first part of proving weak convergence of the path measures. The second part is tightness of the measure family. Tightness will be discussed in the sequel.

\section{Tightness of the conditional process and its projection}

In this first section, we prove that the family of conditional processes is tight, iff the family of its projections onto $L$ is tight.

\subsection{The projection} Let $Y$ be a process on $M$ with continuous paths and $Y^{\eps}$  the process $Y$ conditioned not to leave $L(\eps)$ up to a fixed time $T>0$.  As above, $\pi :L(1)\to L$ denotes the projection. We denote by $S^{\eps}:= \pi Y^{\eps}$ the projected process supported by the path space of $L$. 

\subsection{Eventually basic functions} A continuous function $f\in C(L(1))$ is called {\em eventually basic} if there is some $\delta > 0$ such that $f\vert_{L(\delta)}=f\vert_L\circ\pi$. The subalgebra of eventually basic functions  is denoted by $\BB \subset C(L(1))$ and separates points in $\overline{L(1)}$.

\begin{Lemma}\label{dense} For compact $L$, $\BB\subset C(\overline{L(1)})$ is a dense subset.
\end{Lemma}

\noindent{\bf Proof.} Let $\delta > 0$ and $\phi_1,\phi_2\in C (M)$ be a continuous partition of unity such that
\begin{enumerate}
\item $\phi_1\vert_{L(\delta /2)}=1$ and $\phi_1\vert_{M\setminus L(\delta)}=0$,
\item $\phi_1,\phi_2\geq 0$,
\item $\phi_1(x) + \phi_2 (x) =1$ for all $x\in\overline{L(1)}$.
\end{enumerate}
We consider an arbitrary $f\in C(\overline{L(1)})$. The function $\widehat{f} := (f\vert_L\circ\pi)\phi_1 + f\phi_2$ is eventually basic with
$$
\sup_{x\in L(1)} \vert \widehat{f} - f \vert = \sup_{x\in L(\delta)} \vert (f\vert_L\circ\pi - f)\phi_2\vert \leq \sup_{x\in L(\delta)} \vert f\vert_L\circ\pi - f\vert .
$$
Due to the compactness of $L$, for all $\eps > 0$ we find a $\delta > 0$ such that 
$$
\sup_{x\in L(\delta)} \vert f\vert_L\circ\pi - f\vert < \eps.
$$
That implies the statement.\hfill$\Box$

\subsection{Tightness} From now on, we assume $Y^{\eps}(0) = x\in L$. The conditioned process $Y^{\eps}$ is supported by the path space $\Omega^{\eps}_{0,T}$. $L$ is compact, therefore $\overline{L(1)}$ as well and the {\em compact containment condition} is fulfilled. By \cite{EthKur:86}, 9.1 Theorem, p. 142, together with Lemma \ref{dense}, the family $(Y^{\eps}:\eps > 0)$ is tight, iff for all eventually basic functions $f\in\BB$ the family
$$
\Pi(f) := ( f\circ Y^{\eps}\,:\, \eps > 0)
$$
is tight. \\

\noindent Note that under the assumptions on the tube, the probability of Brownian motion not to leave an $\eps$-tube around $L$ up to time $T>0$ is positive, if the starting point is inside the tube. In particular, for a starting point $x\in L$, this holds for every $\eps$-tube, i.e. $\mu^x (\tau_{\eps}> T) > 0$ for every $\eps > 0$.\\

\begin{Lemma}\label{continuity} For $\eps_0 > 0$, the family $( Y^{\eps}\,:\, \eps_0<\eps \leq 1)$ is tight.
\end{Lemma}

\noindent{\bf Proof} If $\mu^{x}$ denotes the law of $Y$ on $M$ with starting point $x\in L$, the law $\mu_{\eps}^{x}$ of $Y^{\eps}$ with starting point $x$ is given by the Radon - Nikodym density
$$
\Theta_{\eps}= \frac{d\mu_{\eps}^x}{d \mu^x} = \frac{\mathbf{1}_{\OmT}}{\mu^{x} (\OmT)}= \frac{\mathbf{1}_{\OmT}}{\mu^x (\tau_{\eps} > T)},
$$
where $\tau_{\eps}$ denotes the first exit time of the process $Y$ from the tube domain $L(\eps)$. As $\eps$ tends to $\eps_0$, the densities $\Theta_{\eps}$ converge to $\Theta_{\eps_0}$ pointwise. By monotonicity of the first exit time, we have 
$$
0 < \mu^x (\tau_{\eps_0} > T)\leq \mu^x (\tau_{\eps} > T)
$$ 
and hence, by $\Theta_{\eps} \leq  \frac{1}{\mu^x (\tau_{\eps_0} > T)}$
and dominated convergence, we obtain the result.\hfill$\Box$

\begin{Proposition}\label{tightproj} Under the conditions above, the following statements are equivalent:
\begin{enumerate}
\item $( S^{\eps}=\pi Y^{\eps}\,:\, 0<\eps \leq 1)$ is tight,
\item $( Y^{\eps}\,:\, 0<\eps \leq 1)$ is tight.
\end{enumerate}
\end{Proposition}

\noindent{\bf Proof} (a) Let first $( Y^{\eps}\,:\, 0<\eps \leq 1)$ be tight. By the continuous mapping principle, that implies that $( S^{\eps}=\pi Y^{\eps}\,:\, 0<\eps \leq 1)$ is tight as well. (b) Let now $(S^{\eps}=\pi Y^{\eps}\,:\, 0<\eps \leq 1)$ be tight. By the theorem from \cite{EthKur:86} mentioned above, the proposition is proved if we show that this implies that the families $\Pi(f)$ are tight for all $f\in \BB$. Let $\delta > 0$ and  $f\in \BB$ such that there $f\vert_{L(\delta)} = f\vert_L \circ\pi$. By Lemma \ref{continuity}, the set $( Y^{\eps}\,:\, \delta/2 <\eps \leq 1)$ and therefore the set $( f\circ Y^{\eps}\,:\, \delta/2 <\eps \leq 1)$ is tight as well. For $\eps < \delta$, we have
$$
f\circ Y^{\eps} =  f\vert_L\circ\pi Y^{\eps} = f\vert_L\circ S^{\eps}.
$$
By $f\vert_L\in C(L)$, tightness of the family
$$
( f\circ Y^{\eps}\,:\, 0<\eps < \delta) =( f\vert_L\circ S^{\eps}\,:\, 0<\eps \leq \delta)
$$
follows again by the theorem from \cite{EthKur:86} and by our assumption on the tightness of $S^{\eps}$. Thus, by
$$
\Pi (f) = ( f\circ Y^{\eps}\,:\, 0<\eps < \delta ) \cup ( f\circ Y^{\eps}\,:\, \delta/2 <\eps \leq 1)
$$
tightness of $\Pi (f)$ is finally proved.\hfill$\Box$

\section{A sub-Gaussian estimate and tightness of the projections}\label{3}

Denote by $\overline{f}\in C^{\infty}(L(1))$ smooth basic functions $\overline{f} := f\circ\pi$ with $f\in C^{\infty}(L)$. \\

\noindent In the sequel, we consider the pointwise norm of the gradient of functions $g\in C^{\infty}(L(1))$, given by
\begin{equation}\label{gradnorm}
x\mapsto \Vert dg \Vert_{\eps} = \sqrt{\langle dg,dg\rangle_{\eps}} = \sqrt{\star_{\Sa} \, (dg \wedge \#_{\eps} dg)}.
\end{equation}
and the pointwise norm of the gradient of functions $f\in C^{\infty}(L)$, denoted by $\Vert df\Vert_L:=\langle df,df\rangle_L^{1/2}$.

\begin{Lemma}\label{distance} We have
$$
d_L(x,y) = \sup\lbrace \vert f(x) - f(y)\vert\,:\, f\in C^{\infty}(L), \langle df,df\rangle_L\leq 1\rbrace . 
$$
\end{Lemma} 

\noindent{\bf Proof.} See \cite{Dav:87}, proof of Corollary 6, p. 326.\hfill $\Box$\\

\noindent{\bf Remark.} By substituting $cf$ instead of $f$ for an arbitrary $c>0$, the statement of the Lemma also holds with $c d_L(x,y)$ instead of $d_L(x,y)$ and $\langle df,df\rangle_L\leq c^2$ instead of $\langle df,df\rangle_L\leq 1$.\\

\noindent The following proposition, the proof of which will take the entire Section \ref{beweis}, is the key result to prove tightness of the process family. Please note that the usage of the letter $W$ for points on the tube is here to indicate, that we will think of the tube as the neighbourhood $NL(1)$ of the zero section in the normal bundle. For some of the expressions (in particular for $\sigma_{\eps}$) it is more convenient to make use of the vector space structure of the fibers of the normal bundle. 

\begin{Proposition}\label{preparation} Let $t>0$, $\eps > 0$, $h\in C^{\infty}(L)$ with $\sup_{W\in L(1)}\Vert d\overline{h}\Vert_{\eps}\leq 1$ and $K_t^{\eps}$ the kernel given by
\begin{equation*}
e^{-\frac{t}{2}H_{\eps}^0}f (W) = \int_{L(1)} K_t^{\eps} (W,W') f(W') d\mu_{\Sa}(W') .
\end{equation*}
Then we obtain the kernel estimate
\begin{equation}\label{subgaussian_h}
K_t^{\eps}(W,W')\leq Ct^{- \frac{m+3}{2}}\phi_{\eps}(W)\phi_{\eps}(W')e^{-\frac{(\overline{h}(W) - \overline{h}(W'))^2}{4kt}}
\end{equation}
with $k,C>0$.
\end{Proposition}

\noindent Proposition \ref{preparation} implies the following uniform sub-Gaussian kernel estimate:

\begin{Corollary}\label{subgauss} There is some $\eps_0> 0$ such that for all $\eps<\eps_0$, with the notations above, we have 
\begin{equation}\label{subgaussian_d}
K_t^{\eps}(W,W')\leq Ct^{- \frac{m+3}{2}}\phi_{\eps}(W)\phi_{\eps}(W') e^{-\frac{d_L(\pi(W),\pi(W'))^2}{4Bt}}
\end{equation}
with $B,C>0$.
\end{Corollary}

\noindent{\bf Proof.} Let $\overline{h}=h\circ\pi$, $h\in C^{\infty}(L)$ be a basic function. Then, by the expression for the metric on the tube from \cite{NobWit:19}, Proposition 6, we obtain
$$
\langle d\overline{h},d\overline{h}\rangle_{\eps} = \langle dh,dh\rangle_L + \eps \, r_{\eps}(d\overline{h},d\overline{h}) .
$$
Now, by $d\overline{h} = \pi^* dh$, and the fact that the pullback induces a continuous map $\pi^*:C^{\infty }(L)\to C^{\infty}(L(1))$, we obtain
\begin{equation*}
\left\vert \langle dh,dh\rangle_{L} - \langle d\overline{h},d\overline{h}\rangle_{\eps}\right\vert \leq R\eps \,\langle dh,dh\rangle_L
\end{equation*}
with 
$$
R:=\mathrm{sup}_{h\in C^{\infty}(L), h\neq 0} \left( \frac{r_{\eps}(\pi^*dh,\pi^* dh)}{\langle dh, dh\rangle_L}\right) > 0.
$$
That implies for all $0 < \eps < \max\lbrace 1/R,1\rbrace$
\begin{equation}\label{basic_ab}
  (1-R\eps)\,\langle dh,dh\rangle_L \leq  \langle d\overline{h},d\overline{h}\rangle_{\eps}  \leq (1+ R\eps) \,\langle dh,dh\rangle_L
\end{equation}
Now define 
$$
\mathbb{M}_{\eps} = \lbrace \overline{h}= h\circ\pi\in C^{\infty}(L(1))\,:\, \langle dh,dh \rangle_L \leq \frac{1}{1+R\eps}\rbrace .
$$
Thus, by the remark following Lemma \ref{distance}, we have for fixed $x,y\in L$ and all $\eta > 0$ some function $h\in C^{\infty}(L)$ with $\langle dh,dh\rangle_L\leq d_L(x,y)/(1+\eps R)$ such that
$$
\vert h(x) - h(y)\vert \geq \frac{1}{\sqrt{1+\eps R}} d_L(x,y) -\eta\geq \frac{1}{\sqrt{1+R}}d_L(x,y)-\eta.
$$ 
Then, by (\ref{basic_ab}), $\Vert d\overline{h} \Vert_{\eps} \leq 1$, we may apply Proposition \ref{preparation}, and thus, letting finally $\eta$ tend to zero, there is another constant $B>0$ such that
\begin{eqnarray*}
K_t^{\eps}(W,W')&\leq& Ct^{- \frac{m+3}{2}}\phi_{\eps}(W)\phi_{\eps}(W') e^{-\frac{(h(y) - h(x))^2}{4kt}}\\
&\leq& Ct^{- \frac{m+3}{2}}\phi_{\eps}(W)\phi_{\eps}(W') e^{-\frac{d_L(\pi(W),\pi(W'))^2}{4Bt}}.
\end{eqnarray*}
with $\pi(W)=x$, $\pi (W')=y$.
\hfill $\Box$\\

\noindent To prove tightness of the family of conditional measures, we will use the following well-known sufficient criterion.\\

\begin{Proposition}\label{Kolmogorov} The family of processes $(\pi Y^{\eps}\,:\,0 < {\eps} \leq 1)$ is tight, if for all $0\leq s<t\leq T$, there is some $M>0$ such that the expected distance can be estimated by
$$
E^x\left(d_L(\pi  Y_s^{\eps},\pi Y_t^{\eps})^{2M}\right) \leq K \vert t-s\vert^{1+L}
$$ 
for $K,L > 0$ independent of $\eps > 0$.
\end{Proposition}

\noindent{\bf Proof.}\cite{Kal:97}, Corollary 14.9, p. 261.\hfill $\Box$\\

\noindent The ultracontractivity estimate for the kernel is now sufficient to verify the assumptions of Proposition \ref{Kolmogorov} for the projected processes. That implies:

\begin{Theorem}\label{MainB} The family of processes $(Y^{\eps} \,:\,0 < {\eps} \leq 1)$ is tight.
\end{Theorem}

\noindent{\bf Proof.} By Proposition \ref{tightproj}, it suffices to show tightness of the projections. We will apply Proposition \ref{Kolmogorov}. Assume $0\leq s < t\leq 1$ and $Y^{\eps}(0)=x\in L$ the starting point. The expected distance is given by\\
\begin{eqnarray*}
& & E^x \left\lbrack d( \pi Y^{\eps}_t,\pi Y^{\eps}_s)\right\rbrack^{2M}
=\frac{1}{Z_{\eps}(x)}\int_{L(\eps)} d\mu(W) Q_{\eps}(0,x;s,W) ... \\
&&... \int_{L(\eps)} d\mu(W') Q_{\eps}(s,W;t,W') d(\pi W,\pi W')^{2M} ...\\
&& ...\int_{L(\eps)} d\mu(W'') Q_{\eps}(t,W';1,W''),
\end{eqnarray*}
where $Z_{\eps}(x) := \int_{L(\eps)} d\mu(Z) Q_{\eps}(0,x;1,Z)$.
After the transformation to the $1$-tube (note that $x\in L$), we obtain 
\begin{eqnarray*}
& & E^x \left\lbrack d( \pi x^{\eps}_t,\pi x^{\eps}_s)\right\rbrack^{2M} \\
&=&\frac{1}{\Sigma_{\eps}Z(x,\eps)}\Sigma_{\eps}\int_{L(\eps)} d\mu(W) \Sigma_{\eps}^{-1}Q^{\eps}(0,x;s,W) \Sigma_{\eps}... \\
&&... \int_{L(\eps)} d\mu(W') \Sigma_{\eps}^{-1}Q_{\eps}(s,W;t,W')\Sigma_{\eps}\Sigma_{\eps}^{-1} d(\pi W,\pi W')^{2M} ...\\
&& ...\Sigma_{\eps}\int_{L(\eps)} d\mu(W'') \Sigma_{\eps}^{-1}Q_{\eps}(t,W';1,W'')\Sigma_{\eps}\Sigma_{\eps}^{-1}1\\
&=&\frac{1}{Z(x,\eps)}\int_{L(1)} d\mu_{\Sa}(W) K^{\eps}_s(x,W) ... \\
&&...\int_{L(1)} d\mu_{\Sa}(W') K^{\eps}_{t-s}(W,W')\Sigma_{\eps}^{-1} d(\pi W,\pi W')^{2M} ...\\
&& ...\Sigma_{\eps}\int_{L(1)} d\mu_{\Sa}(W'') K^{\eps}_{1-t}(W',W'')\Sigma_{\eps}^{-1}1(W'')\\
&=&\frac{1}{Z(x,\eps)}\int_{L(1)} d\mu_{\Sa}(W) K^{\eps}_s(x,W) ... \\
&&...\int_{L(1)} d\mu_{\Sa}(W') K^{\eps}_{t-s}(W,W')\,d(\pi W,\pi W')^{2M} ...\\
&& ...\int_{L(1)} d\mu_{\Sa}(W'') K^{\eps}_{1-t}(W',W'')\Sigma_{\eps}^{-1}1(W''),
\end{eqnarray*}
where $Z(x,\eps) = \int_{L(1)} d\mu_{\Sa}(V) K^{\eps}_1(x,V)\Sigma_{\eps}^{-1}1(V)$. Now, by Corollary \ref{subgauss}, we obtain
\begin{eqnarray*}
& & \int_{L(1)} d\mu_{\Sa}(W') K^{\eps}_{t-s}(W,W')\Sigma_{\eps}^{-1} d(\pi W,\pi W')^{2M}...\\
& &...\Sigma_{\eps}\int_{L(1)} d\mu_{\Sa}(W'') K^{\eps}_{1-t}(W',W'')\sqrt{\rho(\eps W'')}\\
&=& \int_{L(1)} d\mu_{\Sa}(W') K^{\eps}_{t-s}(W,W')d(\pi W,\pi W')^{2M}... \\
& & ... \int_{L(1)} d\mu_{\Sa}(W'') K^{\eps}_{1-t}(W',W'')\sqrt{\rho(\eps W'')}\\
&\leq& \int_{L(1)} d\mu_{\Sa}(W') \phi_{\eps}(W') e^{-\frac{d_L(\pi(W),\pi(W'))^2}{4B(t-s)}}\left(\frac{d(\pi W,\pi W')^2}{4B(t-s)}\right)^{M}...\\
& & ...  \int_{L(1)} d\mu_{\Sa}(W'') K^{\eps}_{1-t}(W',W'')\sqrt{\rho(\eps W'')}...\\
& & ...  K(4B(t-s))^{M- \frac{m+3}{2}} \phi_{\eps}(W)\\
&\leq& N \,(t-s)^{M- \frac{m+3}{2}} \phi_{\eps}(W)\int_{L(1)} d\mu_{\Sa}(W') \phi_{\eps}(W')...\\
& & ... \int_{L(1)} d\mu_{\Sa}(W'') K^{\eps}_{1-t}(W',W'')\sqrt{\rho(\eps W'')}\\
&=& N' \,(t-s)^{M- \frac{m+3}{2}}\, \phi_{\eps}(W)\, \langle \phi_{\eps},e^{-\frac{1-t}{2}H_{\eps}^0} \sigma_{\eps}^{-1\,*}\sqrt{\rho}\rangle,
\end{eqnarray*}
by $e^{-x}x^M \leq (M/e)^M$, $N':= K (4BM/e)^M$ and $\phi_{\eps}\geq 0$. Now, $\phi_{\eps}\to \phi_0$ (see Corollary \ref{conver} below) and $\sigma_{\eps}^{-1\,*}\sqrt{\rho}\to 1$ uniformly and in $L^2(L(1),\mu_{\Sa})$. Thus, we obtain by Theorem \ref{Main} 
$$
\langle \phi_{\eps}, e^{-\frac{1-t}{2}H_{\eps}^0} \sigma_{\eps}^{-1\,*}\sqrt{\rho}\rangle\to \langle \phi_{0}, \phi_{0}e^{-\frac{1-t}{2}\Delta_L}1\rangle =1, \,\,\,\, \eps\to 0,
$$
and  
$$
Z(x,\eps) = \left(e^{-\frac{1-t}{2}H_{\eps}^0} \sigma_{\eps}^{-1\,*}\sqrt{\rho}\right)(x) \to \phi_0(x), \,\,\,\,\eps\to 0 .
$$
That implies by introducing a new constant $N>0$
 \begin{eqnarray*}
& & \int_{L(1)} d\mu_{\Sa}(W') K^{\eps}_{t-s}(W,W')\Sigma_{\eps}^{-1} d(\pi W,\pi W')^{2M}...\\
& &...\times\Sigma_{\eps}\int_{L(1)} d\mu_{\Sa}(W'') K^{\eps}_{1-t}(W',W'')\sqrt{\rho(\eps W'')}\\
&\leq& N \,(t-s)^{M- \frac{m+3}{2}}\, \phi_{\eps}(W)
\end{eqnarray*}
for all $0<\eps < \eps_0$. Thus, we obtain finally for the expected distance
\begin{eqnarray*}
& & E^x \left\lbrack d( \pi X^{\eps}_t,\pi X^{\eps}_s)\right\rbrack^{2M}\\
&\leq& \frac{N\,\int_{L(1)} d\mu_{\Sa}(W) K^{\eps}_s(x,W) \phi_{\eps}(W)}{Z(x,\eps)} \,(t-s)^{M- \frac{m+3}{2}}\\
&=& N\,(t-s)^{M- \frac{m+3}{2}}. 
\end{eqnarray*}
Taking $M>\frac{m+5}{2}$ yields the statement. \hfill $\Box$

\begin{Corollary} Let $x\in L$ be an arbitrary starting point on the submanifold and $(\mu_{\eps}^x\,:\,\eps > 0)$ the family of path measures of  
Brownian motion on $M$ conditioned to stay in a tubular $\eps$-neighbourhood around $L$ up to some fixed time horizon $T>0$, i.e.
$$
\mu^x_{\eps}(d\om) = \frac{\W^x(d\om\cap \Omega^{\eps}_{0,T})}{\W^x(\Omega^{\eps}_{0,T})}.
$$
Then $\mu_{\eps}^x$ converges weakly to the path measure
$$
\mu^x_{0}(d\om) = \frac{\exp\left(-\frac{1}{2}\int_0^{T} U(\om(s)) ds\right)\mathbb{W}_L(d\om)}{\int_{\Omega} \exp\left(-\frac{1}{2}\int_0^{T} U(\om'(s)) ds\right)\mathbb{W}_L(d\om')}$$
supported by the path space of $L$. 
\end{Corollary}

\noindent{\bf Proof.} The potential $U\in C^{\infty} (L(r))$ from (\ref{potential}) above is bounded. Therefore, the density
$$
\rho (\omega) := \exp\left(\frac{1}{2}\int_0^T \varphi U(\omega (s))ds\right)\in C(\Omega ,\R),
$$
is continuous and bounded below away from zero, i.e. there are real numbers $0< c <C < \infty$ such that $c < \rho (\omega) < C$. Thus, $\rho^{-1}\in C(\Omega,\R)$ has the same properties. Weak convergence of the measures $\nu_{\eps}^x$ to the Wiener measure $\W_L^x$ on $L$ implies for all bounded continuous $\tau \in C(\Omega,\R)$: 
$$
\int_{\Omega} \tau\rho^{-1}(\omega)\nu^x_{\eps}(d\omega) \to \int_{\Omega} \tau\rho^{-1}(\omega)\W_L^x(d\omega),
$$
and also 
$$
\int_{\Omega} \rho^{-1}(\omega)\nu^x_{\eps}(d\omega) \to \int_{\Omega} \rho^{-1}(\omega)\W_L^x(d\omega),
$$
where 
$$
\int_{\Omega} \rho^{-1}(\omega)\nu^x_{\eps}(d\omega) \geq C^{-1} > 0,
$$ 
for all $\eps > 0$. Thus, for all $\tau \in C(\Omega,\R)$, we have
$$
\frac{\int_{\Omega} \tau\rho^{-1}(\omega)\nu^x_{\eps}(d\omega)}{\int_{\Omega} \rho^{-1}(\omega)\nu^x_{\eps}(d\omega)}\to \frac{\int_{\Omega} \tau\rho^{-1}(\omega)\W_L^x(d\omega)}{\int_{\Omega} \rho^{-1}(\omega)\W_L^x(d\omega)}
$$
which implies weak convergence of conditional Brownian motion by
$$
\frac{\int_{\Omega} \tau\rho^{-1}(\omega)\nu^x_{\eps}(d\omega)}{\int_{\Omega} \rho^{-1}(\omega)\nu^x_{\eps}(d\omega)}= \int_{\Omega} \tau(\omega)\mu^x_{\eps}(d\omega).
$$
\hfill $\Box$

\section{Logarithmic Sobolev inequalities and the proof of Proposition \ref{preparation}}\label{beweis}

\subsection{Properties of the ground states}\label{beweis}

Let $q_{\eps}$ be the quadratic form (\ref{qform}) and
$$
q_{\eps}^0 (f) := q_{\eps}(f) - \frac{\lambda_0}{\eps^2}\Vert f\Vert^2.
$$
We consider the mappings $\phi_{\eps,\alpha,0}:\bsob^1(L(1),\mu_{\Sa})\to\R$ from \cite{NobWit:19}, 2.3 with $\alpha > 0$ and
$$
f\mapsto \frac{1}{2} \left(q_{\eps}^{0} (f) +\alpha\Vert f\Vert^2\right),
$$
and $\phi_{\alpha,0}:\bsob^1 (L(1),\mu_{\Sa})\to \R$ with
$$
f\mapsto\left\lbrace\begin{array}{ll} \frac{1}{2} \int_L \left(\langle df_b,df_b\rangle + \alpha f_b^2\right) d\mu_L & ,f\in E_0\\
\infty & ,\mathrm{else}\end{array}\right. ,
$$
where $f\in E_0$ can be uniquely written $f=(f_b\circ\pi)\phi_0$ with $f_b\in C^{\infty}(L)$ (cf. Section \ref{setup}.c). We have

\begin{Proposition}\label{propsoffuncts} For $\alpha > \lambda_0 + 1$ there is some $\eps_0 > 0$ such that for all $\eps < \eps_0$, we have an $A>0$ with
\begin{enumerate}
\item $\phi_{\eps,\alpha,0}(f) \geq \frac{1}{2A}\Vert f\Vert^2_{\bsob^1(L(1),\mu_{\Sa})}$, i.e. the family is equi - coercive,
\item $\phi_{\eps,\alpha,0}$ epi-converges to $\phi_{\alpha,0}$ with respect to the weak topology on $\bsob^1(L(1),\mu_{\Sa})$, as $\eps$ tends to zero. In particular, we have for any sequence $\eps_n$ with $\eps_n\to 0$
\begin{enumerate}
\item $\liminf_{n} \phi_{\eps_n,\alpha,0}(f_n)\geq \phi_{\alpha,0}(f)$ for all sequences $(f_n)_{n\geq 1}$ with $f_n\rightharpoonup f$ weakly in $\bsob^1 (L(1),\mu_{\Sa})$.
\item $\lim_{n} \phi_{\eps_n,\alpha,0}(f)= \phi_{\alpha,0}(f)$, i.e. the functions converge even pointwise.
\end{enumerate}
\item There is some $K>0$ such that $\phi_{\eps,\alpha,0}(\phi_0) \leq \alpha + K\eps$ for all $\eps < \eps_0$. 
\end{enumerate} 
\end{Proposition}

\noindent{\bf Proof.} {\rm (i)} \cite{NobWit:19}, Proposition 2. {\rm (ii)} \cite{NobWit:19}, Proposition 3.\hfill $\Box$\\

\noindent The operators $H_{\eps}$ are bounded below with semi-simple spectrum for all $\eps > 0$. The eigenspace associated to the smallest eigenvalue $\lambda_{\eps}$ is one-dimensional and there is a unique, non-negative eigenfunction $\phi_{\eps}\geq 0$ with $\Vert \phi_{\eps}\Vert = 1$. We are now going to establish uniform upper and lower estimates for the eigenfunctions in terms of the distance to the boundary of $L(1)$ by extending well known bounds for ground states of Dirichlet operators to the parametric situation.\\

\begin{Proposition}\label{asymptotics} We have 
\begin{enumerate}
\item $\lim_{\eps\to 0}\left( \lambda_{\eps}- \frac{\lambda_0}{\eps^2}\right) = 0$,
\item $\lim_{\eps\to 0}\phi_{\eps} = \phi_0$ in $L^2(L(1),\mu_{\Sa})$.
\end{enumerate}
\end{Proposition}

\noindent{\bf Proof.} By the Rayleigh - Ritz principle 
$$
\lambda_{\eps} = \min_{f\in S} \left\lbrack q_{\eps}^0 (f)\right\rbrack ,
$$
where $S=\lbrace \varphi\in\bsob^1(L(1),\mu_{\Sa})\,:\, \Vert \varphi\Vert = 1\rbrace$ and $\Vert - \Vert$ denotes the norm on $L^2(L(1),\mu_{\Sa})$. 

\noindent The mappings $\phi_{\eps,\alpha,0}$ and $\phi_{\alpha,0}$ remain lower semi - continuous when restricted to $S$ when $S$ is equipped with the (weak) relative topology. By Proposition \ref{propsoffuncts} (3) and (1), we obtain that the set of minimizers $M:= \cup_{\eps\leq \eps_0} \mathrm{argmin} \left(\phi_{\eps,\alpha,0}\vert_S\right) \subset S$ is norm - bounded in $\bsob^1(L(1),\mu_{\Sa})$ and therefore contained in a ball $\mathbb{B} = \mathbb{B}(0,r)\subset \bsob^1(L(1),\mu_{\Sa})$ of radius $r>0$. Since $\bsob^1(L(1),\mu_{\Sa})\subset L^2 (L(1),\mu_{\Sa})$ and the embedding is compact, \cite{Dal:93}, Proposition 8.10, p. 93, and Example 8.9, p. 92, respectively, the weak topology on bounded subsets f $\bsob^1(L(1),\mu_{\Sa})$ is induced by the norm on $L^2(L(1),\mu_{\Sa})$ and therefore, the restrictions $\phi_{\eps,\alpha,0}\vert_{S\cap\mathbb{B}}$ are still equi - coercive  by Proposition \ref{propsoffuncts} (1) and epi - converge to $\phi_{\alpha,0}\vert_{S\cap\mathbb{B}}$ by Proposition \ref{propsoffuncts} (2a), (2b). 

\noindent {\rm (i)} Thus, by \cite{Dal:93}, Theorem 7.8, p. 71, for a sequence $(\eps_n)_{n\geq 1}$ with $\eps_n\to 0$ as $n\to\infty$, 
$$
\alpha = \min_{f\in S} \phi_{\alpha,0}(f) = \lim_{n\to\infty}\,\inf_{f\in S} \phi_{\eps_n,\alpha,0}(f) = \lim_{n\to\infty}\left(\lambda_{\eps_n} -\frac{\lambda_0}{\eps_n^2}+ \alpha\right)
$$
and that implies the first statement. {\rm (ii)} $S\cap \mathbb{B}$ is bounded in $\bsob^1(L(1),\mu_{\Sa})$, hence compact in $L^2(L(1),\mu_{\Sa})$. 
We choose a sequence of non-negative ground states $\phi_{\eps_n}\geq 0$ of the operators $H_{\eps_n}$. By compactness, every subsequence of $(\phi_{\eps_n})_{n\geq 1}$ contains a convergent subsequence. By Proposition \ref{propsoffuncts} (2), the limit of this final sequence is a minimizer $\phi^*$ of $\phi_{\eps,\alpha}\vert_{S}$ and, necessarily, $\phi^*\geq 0$. Hence $\phi^*= \phi_0$, and that implies the second statement. \hfill $\Box$

\begin{Corollary}\label{conver} The ground states converge uniformly in all derivatives, i.e. we have 
$$
\lim_{\eps\to 0} \phi_{\eps} = \phi_{0}
$$
in all Sobolev spaces $\sob^n\cap\bsob^1(L(1),\mu_{\Sa})$.
\end{Corollary}

\noindent{\bf Proof.} By Proposition \ref{asymptotics}, (2), we have $\lim_{\eps\to 0} \phi_{\eps}= \phi_0$ in $L^2(L(1),\mu_{\Sa})$. Hence, by \cite{NobWit:19}, Theorem 1, for $t=1$, and a sequence $(\eps_n)$ tending to zero as above, we have 
$$
\lim_{n\to\infty} e^{-\frac{\lambda_{\eps_n}-\eps_n^{-2}\lambda_0}{2}} \phi_{\eps_n} =  \lim_{n\to\infty}e^{-\frac{1}{2}H_{\eps_n}^0}\phi_{\eps_n} = E_0 e^{-\frac{1}{2}\Delta_L}E_0\phi_0 = \phi_0.
$$
uniformly in all derivatives. \hfill $\Box$\\

\noindent Let now $\delta : L(1)\to\R$, $x\mapsto d_{\Sa}(x,\partial L(1))$ denote the distance to the boundary of the tube. It is well known that the non-negative ground state $\varphi$ for the Dirichlet Laplacian on the flat unit ball $\mathbb{D}\subset\R^n$ satisfies an estimate of the form
\begin{equation*}
a' d(x,\partial\mathbb{D}) \leq \varphi (x) \leq A' d(x,\partial\mathbb{D}),
\end{equation*}
with constants $0< a' < A'$. That implies the existence of constants $0<a < A$ such that
\begin{equation}\label{phinull}
a \delta(x) \leq \phi_0 (x) \leq A \delta (x).
\end{equation}
From the uniform convergence of the ground states, we may now conclude that this estimate actually holds uniformly for small values of $\eps$.\\ 

\begin{Corollary}\label{unifestground} There are constants $C>c>0$, such that
$$
c\delta \leq \phi_{\eps} \leq C \delta 
$$
for all $\eps < \eps_0$.
\end{Corollary}

\noindent{\bf Proof.} Let $0<a<A$ be as in (\ref{phinull}). By Corollary \ref{conver}, we may choose $\eps_0 > 0$ so small that
$$
\sup_{x\in L(1)} \Vert d\phi_{\eps} - d\phi_0\Vert_x \leq \frac{a}{2}
$$ 
for all $\eps < \eps_0$. By 
$$
\phi_{\eps}(x) = \int_{\gamma} d\phi_{\eps}
$$
for an arbitrary smooth curve connecting $x$ with a point $x_0\in\partial L(1)$ and 
$$
\left\vert\phi_{\eps}(x) -\phi_0(x)\right\vert = \left\vert\int_{\gamma} d\phi_{\eps} - d\phi_0\right\vert \leq \frac{a}{2}\mathrm{length}(\gamma),
$$
we obtain $\left\vert\phi_{\eps}(x) -\phi_0(x)\right\vert\leq \frac{a}{2}\delta (x)$ by choosing a shortest connection $\gamma$ of $x_0$ to the boundary. Thus
$$
\frac{a}{2}\delta \leq \phi_{\eps}(x) \leq \left(A +\frac{a}{2}\right)\delta,
$$
and that implies the statement.\hfill $\Box$\\

\noindent Finally, we prove a Hardy - inequality for the Dirichlet Laplacian on the tube with respect to the Sasaki metric.\\

\begin{Proposition}{\bf (Hardy - inequality)}\label{Hardy} For $f\in \bsob^1(L(1),\mu_{\Sa})$, we have
\begin{equation}
\int_{L(1)} df\wedge \star_{\Sa} df \geq \frac{1}{4}\int_{L(1)}\star_{\Sa}\frac{f^2}{\delta^2}
\end{equation}
where $\delta : L(1)\to\R$, $x\mapsto d_{\Sa}(x,\partial L(1))$, denotes the distance to the boundary.
\end{Proposition}

\noindent{\bf Proof.} The Hardy inequality is valid for the Dirichlet Laplacian on the flat unit ball $\mathbb{B}\subset \R^{m-l}$ (cf. for instance \cite{Dav:89}, Lemma 1.5.2, p. 26). The fibers $\pi^{-1}(x)$, $x\in L$, with the metric induced by the Sasaki metric are isometric to $\mathbb{B}$, hence (note that $f\in\bsob^1(L(1),\mu_{\Sa})$ implies $f\vert_{\pi^{-1}(x)}\in \bsob^1(\pi^{-1}(x), \mu_x)$ for $\mu_L$-almost all $x\in L$)
\begin{eqnarray*}
\int_{\pi^{-1}(x)}d\mu_x\langle df, df\rangle_{\pi^{-1}(x)} &=& \int_{\mathbb{B}} \Vert df\Vert^2 dx \\
&\geq&\int_{\mathbb{B}}\frac{f^2 dx}{4(1-\Vert x\Vert)^2}= \frac{1}{4}\int_{\pi^{-1}(x)}d\mu_x\frac{f^2}{\delta_x^2}.
\end{eqnarray*}
By \cite{NobWit:19}, the Sasaki metric on one-forms is given by
$$
\alpha\wedge\star_{\Sa}\beta = \langle J_W^*\alpha, J_W^*\beta\rangle_{N^*L} + \langle \kappa_W\alpha, \kappa_W\beta\rangle_L,
$$ 
where $J_W^*$ is the dual of the flat parallel transport $J_W: N_{\pi (W)}L \to T_WN_{\pi(W)}L$ and $\kappa_W:T^*_WNL\to T^*_{\pi(W)}L$ is the dual of the horizontal lift. That implies for $W\in NL$, $x=\pi(W)$,
\begin{eqnarray*}
\int_{L(1)} df\wedge \star_{\Sa} df &=& \int_{L(1)} \left( \langle J^*_Wdf, J_W^* df\rangle_{N^*L} + \langle \kappa_Wdf, \kappa_W df \rangle_L\right) d\mu_{\Sa}\\
&=& \int_{L(1)} \langle J^*_Wdf, J_W^* df\rangle_{N^*L} d\mu_{\Sa}\\
&=&\int_{L} d\mu_{L}(x)\int_{\pi^{-1}(x)}d\mu_x\langle df, df\rangle_{\pi^{-1}(x)}\\
&\geq& \frac{1}{4}\int_{L} d\mu_{L}(x)\int_{\pi^{-1}(x)}d\mu_x\frac{f^2}{\delta_x^2},
\end{eqnarray*}
with $\delta_x:\pi^{-1}(x)\cap L(1)\to\R, \delta_x (W):= d(W, \partial L(1)\cap \pi^{-1}(x))$. By
$\delta (W) = \delta_{\pi(W)}(W)$, and
$$
\int_{L} d\mu_{L}(x)\int_{\pi^{-1}(x)}d\mu_x\frac{f^2}{\delta_x^2} = \int_{L} d\mu_{L}(x)\int_{\pi^{-1}(x)}d\mu_x\frac{f^2}{\delta^2}= \int_{L(1)} \frac{f^2}{\delta^2}d\mu_{\Sa},
$$
we obtain the statement.\hfill $\Box$\\

\subsection{ A uniform estimate for the $L^2$/$L^{\infty}$-norm}

In the sequel, we always assume $\alpha \geq \lambda_0 + 1$. First of all, we prove that the semigroups generated by $H_0^{\eps} + \alpha$ are ultracontractive and prove a uniform bound for the $L^2$/$L^{\infty}$-norm
$\Vert - \Vert_{2,\infty}$ of the operators.

\begin{Proposition}\label{Sobolev1}
There is some $\eps_0 > 0$ and a constant $N>0$, such that for all $\eps < \eps_0$, $t>0$,
$$
\Vert e^{-\frac{t}{2}\left( H_{\eps}^0 + \alpha \right)}\Vert_{2,\infty} \leq N\, t^{-\frac{m+1}{4}},
$$
where $m=\dim M$.
\end{Proposition} 

\noindent{\bf Proof.} By \cite{NobWit:19}, Corollary 4, $\alpha \geq \lambda_0 + 1$ implies $\lambda_s(\eps) + \alpha \geq 1$ for all eigenvalues $(\lambda_s(\eps))_{s\geq 0}$ of $H_{\eps}^0$. By the spectral theorem, we hence obtain for $f\in L^2(L(1),\mu_{\Sa})$ and $t>0$
$$
e^{-\frac{t}{2}\left( H_{\eps}^0 + \alpha \right)}f = \sum_{s\geq 0} e^{-\frac{t}{2}\left( \lambda_s(\eps) + \alpha \right)}\langle u_s(\eps),f\rangle u_s(\eps), 
$$
where the functions $u_s(\eps)\in L^2(L(1),\mu_{\Sa})$, $s\geq 0$ are smooth and normalized eigenfunctions. By \cite{NobWit:19}, Corollary 10, we have for all $n\in\N$ some $D_n>0$ such that we have the following estimate for the $2n$-Sobolev norm on $\sob^{2n}(L(1),\mu_{\Sa})$:
\begin{eqnarray*}
& & \left\Vert e^{-\frac{t}{2}\left( H_{\eps}^0 + \alpha \right)}f\right\Vert^2_{2n} \\
&\leq& 2D_n^2\left(\Vert e^{-\frac{t}{2}\left( H_{\eps}^0 + \alpha \right)}f\Vert^2 + \Vert \left( H_{\eps}^0 + \alpha \right)^{2n}e^{-\frac{t}{2}\left( H_{\eps}^0 + \alpha \right)}f\Vert^2 \right)\\
&=&  2D_n^2 \sum_{s\geq 0} \left(1 + (\lambda_s(\eps)+\alpha)^{2n}\right)e^{-t(\lambda_s(\eps)+ \alpha)}\left\vert\langle u_s(\eps),f\rangle\right\vert^2 .\end{eqnarray*}
Now, for $x\geq 1$ and $t>0$, we have the estimate
$$
x^{2n} e^{-x} \leq \chi_n^2(t) := \left\lbrace\begin{array}{ll} \left(\frac{2n}{et}\right)^{2n} & ,t\leq 2n \\
e^{-t} & ,t > 2n \\
\end{array}\right.\leq c_n^2 t^{-2n},$$
where $c_n := (2n/e)^{n}$ and $\chi_n:(0,\infty)\to\R$ is continuous and strictly decreasing. Therefore, 
$$
\left\Vert e^{-\frac{t}{2}\left( H_{\eps}^0 + \alpha \right)}f\right\Vert^2_{2n}\leq 2D_n^2 (e^{-t} + \chi_n^2(t)) \Vert f\Vert^2,
$$
or 
\begin{equation}\label{soboab_even}
\left\Vert e^{-\frac{t}{2}\left( H_{\eps}^0 + \alpha \right)}f\right\Vert_{2n}\leq 2D_n \chi_n(t) \Vert f\Vert \leq K'_{2n} t^{-n}\,\Vert f\Vert,
\end{equation}
with $K'_{2n} := 2\sqrt{2} D_nc_n$. By interpolation, we obtain for the $2n+1$-Sobolev norm the estimate
\begin{equation}\label{soboab_odd}
\left\Vert e^{-\frac{t}{2}\left( H_{\eps}^0 + \alpha \right)}f\right\Vert_{2n+1}  \leq K'_{2n+1} t^{-\frac{2n+1}{2}}\,\Vert f\Vert,
\end{equation}
with some suitable $K'_{2n+1} > 0$. \\

\noindent By the Sobolev embedding theorem, we have for $2r>m$ some constant $A_r > 0$ such that
$$
\Vert u\Vert_{\infty} \leq A_r \Vert u\Vert_r .
$$
That implies for $r=\frac{m+1}{2}$, 
\begin{eqnarray*}
\left\Vert e^{-\frac{t}{2}\left( H_{\eps}^0 + \alpha \right)}f\right\Vert_{\infty}&\leq& A_{\frac{m+1}{2}} \left\Vert e^{-\frac{t}{2}\left( H_{\eps}^0 + \alpha \right)}f\right\Vert_{\frac{m+1}{2}} \\
&\leq& K_m t^{-\frac{m+1}{4}}  \Vert f\Vert\end{eqnarray*}
where $K_m := A_r K'_{r} > 0$, and we use either (\ref{soboab_even}) or (\ref{soboab_odd}) depending on whether $m$ is even or odd. \hfill $\Box$\\

\noindent By the symmetry mentioned above, it will suffice, if we consider for all $\eps > 0$ the semigroup generated by $H_{\eps}$, i.e. the operator associated to the (non - renormalized) quadratic form (cf. (\ref{qform}))
\begin{equation*}
q_{\eps} (f) := \int_{L(1)} df\wedge \#_{\eps}df,\,\,\,\,\, \DD = \bsob^1 (L(1),\mu_{\Sa}) .
\end{equation*}
Clearly, the same estimate as in Proposition \ref{Sobolev1} also holds for the semigroup generated by $H_{\eps}$:\\

\begin{Corollary}\label{Sobolev2} Assume $\alpha > \lambda_0+1$. There is some $\eps_0 > 0$ and a constant $N>0$, such that for all $\eps < \eps_0$
$$
\Vert e^{-\frac{t}{2} H_{\eps}}\Vert_{2,\infty} \leq N\, t^{-\frac{m+1}{4}}e^{-\frac{t}{2}\left(\frac{\lambda_0}{\eps^2}-\alpha\right)},
$$
where $m=\dim M$.
\end{Corollary}

\noindent{\bf Proof.} By
$$
N\, t^{-\frac{m+1}{4}}\geq\Vert e^{-\frac{t}{2}\left( H_{\eps}^0 + \alpha \right)}\Vert_{2,\infty} = e^{\frac{t}{2}\left( \frac{\lambda_0}{\eps^2} -\alpha\right)}\Vert e^{-\frac{t}{2} H_{\eps}}\Vert_{2,\infty},
$$
for $\alpha > \lambda_0 + 1$, we obtain the statement. \hfill $\Box$\\

\noindent In order to apply the arguments from \cite{Dav:89}, 2.2, we have to show we are actually dealing with {\em symmetric Markov semigroups}, i.e. we will now check the Beurling - Deny conditions.

\begin{Proposition}\label{Markov_SG} For $\eps_0 > \eps > 0$, the semigroup
$$
t\mapsto e^{-\frac{t}{2} H_{\eps}}
$$
is positivity - preserving and a semigroup of $L^{\infty}$-contractions.
\end{Proposition}

\noindent{\bf Proof.} By \cite{NobWit:19}, Corollary 4, $H_{\eps}^0+\alpha$ is self - adjoint and non - negative, since $\alpha \geq \lambda_0 +1$. That implies, $H_{\eps}$ is self - adjoint as well and non - negative, since $q_{\eps}\geq 0$. Hence, the associated extended quadratic form
$$
f\mapsto \left\lbrace\begin{array}{ll} q_{\eps}(f) & ,f\in\bsob^1(L(1),\mu_{\Sa})\\
\infty & , \mathrm{else}\end{array}\right.
$$ 
is lower semicontinuous on $L^2(L(1),\mu_{\Sa})$.
{\rm (i)} Let $f\in \bsob^1(L(1),\mu_{\Sa})$ and $\phi_t\in C^1(\R)$, $t>0$, be a family of functions such that
\begin{enumerate}
\item $\vert \phi_t'(x)\vert\leq 1$ for all $x\in\R$,
\item $0 < \phi_t(x) - \vert x\vert \leq t$ for all $x\in\R$, i.e. $\phi_t$ converges uniformly from above to $\vert x\vert$, as $t$ tends to zero. 
\end{enumerate}
That implies $\phi_t\circ f\in \bsob^1(L(1),\mu_{\Sa})$ (note that $\mu_{\Sa}(L(1))<\infty$) and
\begin{eqnarray*}
q_{\eps} (\phi_t\circ f) &=& \int_{L(1)} d(\phi_t\circ f)\wedge\#_{\eps} d(\phi_t\circ f) = \int_{L(1)} \vert\phi'_t(f)\vert^2 df\wedge\#_{\eps} df\\
&\leq& \int_{L(1)}  df\wedge\#_{\eps} df = q_{\eps}(f).
\end{eqnarray*}
Again, by $\mu_{\Sa}(L(1))<\infty$, uniform convergence $\phi_t\to\vert x\vert$, $t\to 0$, implies $\Vert \phi_t \circ f - \vert f\vert\Vert\to 0$ in $L^2(L(1),\mu_{\Sa})$ and thus,
$$
q_{\eps} (\vert f\vert )\leq \liminf_{t\to 0} q_{\eps} (\phi_t\circ f) \leq q_{\eps}(f),
$$
which implies $\vert f\vert\in\bsob^1(L(1), \mu_{\Sa})$ and finally positivity of the semigroup by \cite{Dav:89}, Thm. 1.3.2, p. 12.
 {\rm (ii)} Let $f\in \bsob^1(L(1),\mu_{\Sa})$ and $\psi_t\in C^1(\R)$, $t>0$, be a family of functions such that 
 \begin{enumerate}
\item $0\leq \psi'_{t}(x)\leq 1$ for all $x\in\R$,
\item $\sup_{x\in\R}\vert\psi_{t}(x) - 0\vee (x\wedge 1)\vert \leq t$ for all $t >0$, i.e. the family approximates $0\vee (x\wedge 1)$ uniformly. 
\end{enumerate}
Therefore, $f\in \bsob^1(L(1),\mu_{\Sa})$  implies $\psi_{t}\circ f\in \bsob^1(L(1),\mu_{\Sa})$ and, since $\mu_{\Sa}(L(1))<\infty$, we have $\Vert  \psi_{t}\circ f - 0\vee (f\wedge 1)\Vert \to 0$ as $t \to 0$ in $L^2(L(1),\mu_{\Sa})$. As above, $q_{\eps} (\psi_t\circ f) \leq q_{\eps}(f)$, and therefore
$$
q_{\eps} (0\vee (f\wedge 1))\leq \liminf_{t\to 0} q_{\eps} (\phi_t\circ f) \leq q_{\eps}(f).
$$
Hence, $0\vee (f\wedge 1)\in\bsob^1(L(1), \mu_{\Sa})$. Thus, \cite{Dav:89}, Thm. 1.3.3, p. 14 finally implies the contraction property.\hfill $\Box$\\

\noindent Corollary \ref{Sobolev2} and Proposition \ref{Markov_SG} immediately imply an $L^2$-Logarithmic Sobolev inequality for $H_{\eps}$. In the sequel, we will denote by
\begin{equation}\label{entropy}
\EE_p(f) := \int_{L(1)} f^p \log\left(\frac{f}{\Vert f\Vert_p}\right), \,\,\,\, \Vert f\Vert_p^p := \int_{L(1)}\star_{\Sa} f^p\end{equation}
the $p$-entropy term and the $p$-norm contained in the Logarithmic Sobolev inequalities, and we always assume $p\geq 2$.

\begin{Proposition}\label{L2_logsobo} For $0\leq f\in\bsob^1(L(1),\mu_{\Sa})\cap L^{\infty}(L(1),\mu_{\Sa})$, we have $f^2\log f\in L^1(L(1),\mu_{\Sa})$ and
\begin{equation}\label{logsoboA}
\EE_2 (f) \leq \theta \frac{1}{2}q_{\eps}(f) + \beta (\theta) \Vert f\Vert_2^2, 
\end{equation}
for all $\theta >0$, where $\beta (\theta) = c - \frac{m+1}{4} \log (\theta) - \frac{\theta}{2}\left(\frac{\lambda_0}{\eps^2}-\alpha\right)$. The estimate holds uniformly for all $0< \eps < \eps_0$. 
\end{Proposition}

\noindent{\bf Proof.} \cite{Dav:89}, Theorem 2.2.3, p. 64.\hfill $\Box$\\

\noindent As in Section \ref{beweis}, we denote by $q_{\eps}^0$ the renormalized quadratic form
$$
q_{\eps}^0 (f) =\int_{L(1)} df \wedge \#_{\eps} df - \frac{\lambda_0}{\eps^2}\Vert f\Vert_2^2,
$$
where $\lambda_0>0$ is the smallest eigenvalue of the Dirichlet Laplacian on the flat unit ball $\mathbb{B}\subset \R^{m-l}$. $q_0^{\eps}$ is the form associated to the operator $H_{\eps}^0$. Now, a simple rearrangement of (\ref{logsoboA}) yields:\\

\begin{Corollary}\label{L2_Logsobo_B} Under the same hypotheses as in Proposition \ref{L2_logsobo}, we have
\begin{equation}\label{logsoboB}
\EE_2 (f) \leq \theta \frac{1}{2}\left(q_{\eps}^0(f) + \alpha \Vert f\Vert_2^2\right) + \gamma (\theta) \Vert f\Vert_2^2, 
\end{equation}
with $\gamma (\theta) = c - \frac{m+1}{4} \log (\theta)$ for all $\theta >0$ uniformly for all $0 < \eps < \eps_0$.
\end{Corollary} 

\subsection{Weighted $L^2$-Sobolev estimates}

From Corollary \ref{L2_Logsobo_B}, we will infer weighted Sobolev estimates in $L^p$-spaces for general $p\geq 2$. We consider a ground state transform together with a weight of the form $e^{a \widehat{h}}$ with $a\in\R$ and $\widehat{h}\in C^{\infty}(\overline{L(1)})$ to be specified later. In the sequel, we denote by
$$
B_{\eps}^0 (f,g) := \int_{L(1)} df \wedge \#_{\eps} dg - \star_{\Sa}\lambda_{\eps} fg
$$
the bilinear form associated to the quadratic form $q_{\eps}(f) - \lambda_{\eps}\Vert f\Vert^2$, where $\lambda_{\eps}>0$ is the smallest eigenvalue of the operator $H_{\eps}$, i.e. $H_{\eps}\phi_{\eps}=\lambda_{\eps}\phi_{\eps}$ as in Section \ref{beweis} above.\\

\noindent Please note that the chosen renormalizations differ for $q_{\eps}^0(f)$ and $B_{\eps}^0(f,g)$, such that there is an overall difference of
\begin{equation}\label{differreno}
q_{\eps}^0(f) = q_{\eps}(f) -\frac{\lambda_0}{\eps^2}\Vert f\Vert_2^2 = B_{\eps}^0(f,f) + \left( \lambda_{\eps} - \frac{\lambda_0}{\eps^2}\right)\Vert f\Vert_2^2 .
\end{equation} 
By Proposition \ref{asymptotics}, $q_{\eps}^0(f) - B_{\eps}^0(f,f)\to 0$ as $\eps$ tends to zero.\\

\noindent By the ground state transform, the measure underlying the Hilbert space under consideration is changed from $\mu_{\Sa}$ to $\phi_{\eps}^2\mu_{\Sa}$ such that we have to adjust our notation. Namely, we will use the shorthands
\begin{equation}
\EE_{\eps,p}(f) := \int_{L(1)} \star_{\Sa}\phi_{\eps}^2 f^p \log\left(\frac{f}{\Vert f\Vert_p}\right), \,\,\,\Vert f\Vert_{\eps,p}^p :=  \int_{L(1)}\star_{\Sa}  \phi_{\eps}^2f^p 
\end{equation}
for $p$-entropy term and $p$-norm ($p\geq 2$), respectively.\\

\noindent Since we want to apply a ground state transform to the bilinear form $B_{\eps}^0(-,-)$, we consider from now on the domain
\begin{equation}
\CC := \lbrace f\in L^2(L(1),\nu_{\eps})\,:\,\phi_{\eps}f \in \bsob^1 (L(1),\mu_{\Sa})\rbrace \cap L^{\infty}(L(1),\nu_{\eps})
\end{equation}
where $\nu_{\eps}$ is the probability measure on $L(1)$ given by the radon - Nikodym density $\frac{d\nu_{\eps}}{d\mu_{\Sa}} = \phi_{\eps}^2$. Clearly, $L^{\infty}(L(1),\nu_{\eps})=L^{\infty}(L(1),\mu_{\Sa})$ and therefore, since $\phi_{\eps}\geq 0$ is bounded, we have
$$
\CC \subset\lbrace f\in L^2(L(1),\nu_{\eps})\,:\,\phi_{\eps}f \in \bsob^1 (L(1),\mu_{\Sa})\cap L^{\infty}(L(1),\mu_{\Sa})\rbrace .
$$

\begin{Lemma}\label{estimate_A} Let $f,g\in\CC$. Then
\begin{enumerate}[label=(\roman*)]
\item $B_{\eps}^0(\phi_{\eps}f,\phi_{\eps}g) = \int_{L(1)} \phi_{\eps}^2 \,df\wedge\#_{\eps}dg$.
\item For $p\geq 2$, we have
$$ 
B_{\eps}^0 (\phi_{\eps}f^{p/2},\phi_{\eps}f^{p/2}) = \frac{p^2}{4(p-1)} B_{\eps}^0 (\phi_{\eps}f,\phi_{\eps}f^{p-1}) .
$$
\item Let $\widehat{h}\in C^{\infty}(\overline{L(1)})$, $a\in\R$. For $p\geq 2$, we have
\begin{eqnarray*}
&& B_{\eps}^0 (e^{a\widehat{h}} \phi_{\eps}f,e^{-a\widehat{h}} \phi_{\eps}f^{p-1}) - B_{\eps}^0 (\phi_{\eps}f,\phi_{\eps}f^{p-1}) \\
&=& - \int_{L(1)}  \phi_{\eps}^2\left( a^2 f^p d\widehat{h}\wedge\#_{\eps} d\widehat{h} - \frac{a (p-2)}{p}df^p\wedge\#_{\eps}d\widehat{h}\right) .
\end{eqnarray*}
\end{enumerate}
\end{Lemma}

\noindent{\bf Proof.} {\rm (i)} $\phi_{\eps}$ vanishes on the boundary of the tube. Hence, by Stokes' theorem,
\begin{eqnarray*}
& & \int_{L(1)} d(\phi_{\eps}f) \wedge \#_{\eps} d(\phi_{\eps}g) - \star_{\Sa}\lambda_{\eps} \phi_{\eps}^2 gf \\
&=&  \int_{L(1)}  \phi_{\eps}^2 df \wedge \#_{\eps} dg + d(fg\phi_{\eps}\#_{\eps}d\phi_{\eps}) +   \star_{\Sa}\,fg(\phi_{\eps} H_{\eps}\phi_{\eps}-\lambda_{\eps} \phi_{\eps}^2)  \\
&=& \int_{L(1)}  \phi_{\eps}^2 df \wedge \#_{\eps} dg.
\end{eqnarray*}
{\rm (ii)} follows from {\rm (i)} and 
$$
df^{p/2}\wedge \#_{\eps} df^{p/2} = \frac{p^2}{4(p-1)} df \wedge \#_{\eps} df^{p-1} .
$$
{\rm (iii)} follows from {\rm (i)}, $fd(\phi_{\eps}f^{p-1}) = f^{p-1}d(\phi_{\eps}f) +\frac{p-2}{p}\phi_{\eps}df^p$, and the identity $e^{\mp a \widehat{h}} d(e^{\pm a \widehat{h}}f) = df \pm a f d\widehat{h}$. \hfill$\Box$

\begin{Lemma}\label{estimate_B} We have for $f\in\CC$, $f\geq 0$, $\tau > 0$, and $\sup_{W\in L(1)}\Vert d\widehat{h}\Vert_{\eps}\leq 1$
$$
\left\vert \int_{L(1)} \phi_{\eps}^2 df^p\wedge\#_{\eps}d\widehat{h}\right\vert \leq \tau  B_{\eps}^0(\phi_{\eps}f^{p/2},\phi_{\eps}f^{p/2}) + \frac{1}{\tau}\Vert f\Vert_{\eps,p}^p .
$$
\end{Lemma}

\noindent{\bf Proof} 
By the Cauchy-Schwarz inequality and $\Vert d\widehat{h}\Vert_{\eps}\leq 1$
\begin{eqnarray*}
\left\vert\int_{L(1)}  \phi_{\eps}^2df^p\wedge\#_{\eps}d\widehat{h}\right\vert &=& \left\vert\int_{L(1)}  \phi_{\eps}^2 \langle 2f^{p/2}df^{p/2},d\widehat{h}\rangle_{\eps}d\mu_{\Sa}\right\vert \\
&\leq& 2\int_{L(1)}   \phi_{\eps}^2 f^{p/2}\langle df^{p/2},df^{p/2}\rangle^{1/2}_{\eps}d\mu_{\Sa},
\end{eqnarray*}
hence for all $\tau >0$
\begin{eqnarray*}
\left\vert\int_{L(1)}  \phi_{\eps}^2df^p\wedge\#_{\eps}d\widehat{h}\right\vert 
&\leq& \int_{L(1)}  \phi_{\eps}^2 \left(\frac{1}{\tau} f^{p} + \tau \langle df^{p/2},df^{p/2}\rangle_{\eps}\right)d\mu_{\Sa}\\
&=& \int_{L(1)}  \phi_{\eps}^2\left(\frac{1}{\tau}\star_{\Sa} f^{p} + \tau df^{p/2}\wedge\#_{\eps}df^{p/2}\right).
\end{eqnarray*}
By Lemma \ref{estimate_A} (i), that implies the statement. \hfill $\Box$\\

\noindent As a corollary, we obtain the following estimate of the weighted transformed bilinear form, which will be one part of the final estimate in the next section.

\begin{Corollary} \label{estimate_C}
We have for $f\in\CC$, $f\geq 0$, $\sup_{W\in L(1)}\Vert d\widehat{h}\Vert_{\eps}\leq 1$ and $p\geq 2$
$$
B_{\eps}^0 (e^{a \widehat{h}}\phi_{\eps}f,e^{-a \widehat{h}}\phi_{\eps}f^{p-1}) + \frac{a^2 p}{2}\Vert f\Vert_{\eps,p}^p \geq \frac{2}{p}B_{\eps}^0(\phi_{\eps}f^{p/2},\phi_{\eps}f^{p/2}) .$$
\end{Corollary}

\noindent{\bf Proof.} By $\Vert d\widehat{h}\Vert_{\eps}\leq 1$, Lemma \ref{estimate_A}, {\rm (iii)}, and Lemma \ref{estimate_B}
\begin{eqnarray*}
&  & B_{\eps}^0 (e^{a \widehat{h}}\phi_{\eps}f,e^{-a \widehat{h}}\phi_{\eps}f^{p-1}) + \rho^2 \Vert f\Vert_{\eps,p}^p \\
&\geq & B_{\eps}^0 (e^{a \widehat{h}}\phi_{\eps}f,e^{-a \widehat{h}}\phi_{\eps}f^{p-1}) + \rho^2 \int_{L(1)}\phi_{\eps}^2 f^p d\widehat{h}\wedge\#_{\eps} d\widehat{h} \\
&=& B_{\eps}^0 (\phi_{\eps}f,\phi_{\eps}f^{p-1}) + \frac{a (p-2)}{p}\int_{L(1)}  \phi_{\eps}^2df^p \wedge\#_{\eps}d\widehat{h} \\
&\geq & B_{\eps}^0 (\phi_{\eps}f,\phi_{\eps}f^{p-1}) - \frac{\vert a\vert\, (p-2)}{p}\int_{L(1)}  \phi_{\eps}^2df^p \wedge\#_{\eps}d\widehat{h} \\
&\geq& B_{\eps}^0 (\phi_{\eps}f,\phi_{\eps}f^{p-1}) - \frac{\vert a\vert\, (p-2)\tau}{p}  B_{\eps}^0(\phi_{\eps}f^{p/2},\phi_{\eps}f^{p/2}) \\
& & - \frac{\vert a\vert\, (p-2)}{p\tau} \Vert f\Vert_{\eps,p}^p.
\end{eqnarray*}
Hence, by Lemma \ref{estimate_A}, {\rm (ii)}
\begin{eqnarray*}
&  & B_{\eps}^0 (e^{a \widehat{h}}\phi_{\eps}f,e^{-a \widehat{h}}\phi_{\eps}f^{p-1}) + \left(a^2 + \frac{\vert a\vert\, (p-2)}{p\tau}\right)\Vert f\Vert_{\eps,p}^p\\
&\geq& \left( \frac{4(p-1)} {p^2} - \frac{\vert a\vert\, (p-2)\tau}{p} \right) B_{\eps}^0(\phi_{\eps}f^{p/2},\phi_{\eps}f^{p/2}) .
\end{eqnarray*}
Furthermore, by $p\geq 2$, the choice $\tau = 2/(p\,\vert a\vert)> 0$ yields
$$
B_{\eps}^0 (e^{a \widehat{h}}\phi_{\eps}f,e^{-a \widehat{h}}\phi_{\eps}f^{p-1}) + \frac{a^2 p}{2}\Vert f\Vert_{\eps,p}^p \geq \frac{2}{p}B_{\eps}^0(\phi_{\eps}f^{p/2},\phi_{\eps}f^{p/2}) .$$
\hfill $\Box$

\noindent The Rosen - Lemma is the last preparation for the weighted Log-Sobolev estimate in the next section.

\begin{Proposition}{\bf (Rosen - Lemma)}\label{Rosen} Let $f\in\CC$, $f\geq 0$ and $p\geq 2$. For all $\tau > 0$, we have 
$$
-\int_{L(1)} \star_{\Sa}\phi_{\eps}^2 f^p \log \phi_{\eps} \leq \tau B_{\eps}^0 (\phi_{\eps} f^{p/2},\phi_{\eps} f^{p/2}) +\nu(\tau)\Vert f\Vert_{\eps,p}^p,$$
where $\nu(\tau) := k_1 + k_2\tau - \frac{1}{2}\log (\tau)$.
\end{Proposition}

\noindent{\bf Proof.} For all $\tau' > 0$, we have $-x\log x \leq g(\tau')x + \tau'$ for all $x > 0$, where $g(\tau') :=  - \log (\tau' e) $ for $\tau' > 0$. Hence 
\begin{eqnarray*}
-\int_{L(1)} \star_{\Sa}\phi_{\eps}^2 f^2 \log \phi_{\eps} &=& -\frac{1}{2}\int_{L(1)} \star_{\Sa}\phi_{\eps}^2 f^2 \log \phi_{\eps}^2 \\
&\leq& \frac{1}{2}g(\tau')\int_{L(1)}\star_{\Sa}\phi_{\eps}^2 f^2 +\frac{\tau'}{2}\int_{L(1)} \star_{\Sa} f^2 \\
&=& \frac{1}{2}g(\tau') \int_{L(1)}\star_{\Sa}\phi_{\eps}^2 f^2 +\frac{\tau'}{2}\int_{L(1)} \star_{\Sa} \phi_{\eps}^2 \frac{f^2}{\phi_{\eps}^2}.
\end{eqnarray*}
By $c\delta \leq \phi_{\eps}\leq C\delta$ from Corollary \ref{unifestground} and the Hardy - inequality from Proposition \ref{Hardy}, we have
$$
\int_{L(1)} \star_{\Sa} \phi_{\eps}^2 \frac{f^2}{\phi_{\eps}^2}\leq \frac{1}{c^2} \int_{L(1)}\star_{\Sa}\frac{\phi_{\eps}^2f^2}{\delta^2}\leq \frac{4}{c^2} \int_{L(1)}d(\phi_{\eps}f)\wedge\star_{\Sa}d(\phi_{\eps}f).
$$
Now by \cite{NobWit:19}, Proposition 2, there is some $A >0$ such that for $\alpha > \lambda_0$
$$
q_{\eps}^0(f) + \alpha \Vert f\Vert_2^2 \geq \frac{1}{2A}\Vert f\Vert_{\bsob^1(L(1),\mu_{\Sa})}.
$$
Hence, by the definition of the $\sob^1$-norm, 
\begin{eqnarray*}
& & \int_{L(1)}d(\phi_{\eps}f)\wedge\star_{\Sa}d(\phi_{\eps}f) \leq \Vert \phi_{\eps}f \Vert_{\bsob^1(L(1),\mu_{\Sa})}^2\\
&\leq& 2A\left(B^0_{\eps}(\phi_{\eps}f,\phi_{\eps}f) + \left(\lambda_{\eps}-\frac{\lambda_0}{\eps^2} +\alpha \right) \Vert f \Vert^2_{\eps,2}\right). 
\end{eqnarray*}
Let now $K':= 8Ac^{-2}>0$ and $\tau = K'\tau'>0$. Then
\begin{eqnarray*}
& & -\int_{L(1)} \star_{\Sa}\phi_{\eps}^2 f^2 \log \phi_{\eps}\leq \frac{1}{2}g(\tau') \Vert f\Vert_{\eps,2}^2 +\frac{\tau'}{2}\int_{L(1)} \star_{\Sa} \phi_{\eps}^2 \frac{f^2}{\phi_{\eps}^2} \\
&\leq&\frac{1}{2}g(\tau') \Vert f\Vert_{\eps,2}^2 +\frac{K'\tau'}{2} \left(B^0_{\eps}(\phi_{\eps}f,\phi_{\eps}f) + \left(\lambda_{\eps}-\frac{\lambda_0}{\eps^2} +\alpha \right) \Vert f \Vert^2_{\eps,2}\right)\\
&=&\frac{\tau}{2} B^0_{\eps}(\phi_{\eps}f,\phi_{\eps}f) + \frac{1}{2}\left( \tau \left(\lambda_{\eps}-\frac{\lambda_0}{\eps^2} +\alpha \right) + g(\tau/K')\right)\Vert f\Vert_{\eps,2}^2 .
\end{eqnarray*}
Now, by Proposition \ref{asymptotics}, there ist some $k_2>0$ with $2^{-1}(\lambda_{\eps}- \lambda_0\eps^{-2} +\alpha) \leq k_2$ uniformly in $\eps < \eps_0$. With $k_1:= -\frac{1}{2}\log(e/K')$, we obtain the statement for $p=2$. The case $p>2$ follows since $f\in\CC$, $f\geq 0$, implies $f^{p/2}\in\CC$. \hfill $\Box$

\subsection{Weighted $L^p$-Sobolev estimates} We will now derive logarithmic Sobolev inequalities for functions $f\in\CC$ with $f\geq 0$. All vectors of the form $e^{-t\HH_{\eps}}u$ with $t>0$ and $u\in L^2(L(1),\nu_{\eps})$ are contained in $\CC$, and to calculate the differential $\frac{d}{ds}\Vert e^{-s\HH_{\eps}}u\Vert_{p(s)}^{p(s)}$, which is essential for the estimate of the semigroup norm, it is therefore sufficient to consider functions $f\in\CC$. Furthermore, since the semigroup preserves positivity, we may even restrict the estimates to non-negative $f\geq 0$. 

\begin{Theorem}\label{LP_logsobo} For $p\geq 2$, $\zeta > 0$, $f\in\CC$ and $\widehat{h}\in C^{\infty}(\overline{L(1)})$ with $\sup_{W\in L(1)}\Vert d\widehat{h}\Vert_{\eps}\leq 1$, we have the estimate
\begin{equation*}
\EE_{\eps,p}(f) \leq \frac{\zeta}{2}B_{\eps}^0 (e^{a \widehat{h}}\phi_{\eps}f,e^{-a \widehat{h}}\phi_{\eps}f^{p-1})  + \gamma (\zeta, p)\, \Vert f\Vert_{\eps,p}^p,\end{equation*}
where
$$
\gamma (\zeta,p) = \frac{\zeta a^2p }{2} +\frac{2}{p}\left(A +  \frac{B\zeta}{p}-\frac{m+3}{4}\log(\zeta/p)\right)
$$
with constants $A,B > 0$.
\end{Theorem}

\noindent{\bf Proof.} The $p$-entropy is given by
\begin{eqnarray*}
\EE_{\eps,p}(f) &=& \frac{2}{p}\left( \EE_2 (\phi_{\eps}f^{p/2}) - \int_{L(1)}\star_{\Sa} \phi_{\eps}^2 f^p \log \phi_{\eps}\right).
\end{eqnarray*}
On the one hand, by Corollary \ref{L2_Logsobo_B} and (\ref{differreno}), we have
$$
\EE_2 (\phi_{\eps}f^{p/2})\leq \frac{\theta}{2} B_{\eps}^0(\phi_{\eps}f^{p/2},\phi_{\eps}f^{p/2})+ \eta (\theta)\Vert f\Vert_{\eps,p}^p
$$
with 
$$
\eta (\theta) = \left(\frac{\theta}{2}\left(\lambda_{\eps}-\frac{\lambda_0}{\eps^2} +\alpha \right)+\gamma(\theta)\right)\leq c_1 + c_2\theta -\frac{m+1}{4}\log(\theta).
$$
On the other hand, by Proposition \ref{Rosen}
$$
- \int_{L(1)}\star_{\Sa} \phi_{\eps}^2 f^p \log \phi_{\eps} \leq \frac{\tau}{2} B_{\eps}^0 (\phi_{\eps}f^{p/2}, \phi_{\eps}f^{p/2}) + \nu (\tau) \Vert f\Vert_{\eps,p}^p,
$$
with $\nu (\tau)=k_1 + k_2\tau -\frac{1}{2}\log (\tau)$. Hence, for $\tau,\theta > 0$,
\begin{eqnarray*}
& & \EE_{\eps,p}(f) - \int_{L(1)}\star_{\Sa} \phi_{\eps}^2 f^p \log \phi_{\eps} \\
&\leq&  \frac{\theta+\tau}{2} B_{\eps}^0(\phi_{\eps}f^{p/2},\phi_{\eps}f^{p/2})+ (\eta (\theta)+\nu(\tau))\Vert f\Vert_{\eps,p}^p,
\end{eqnarray*}
and that implies by Corollary \ref{estimate_C}
\begin{eqnarray*}
\EE_{\eps,p}(f) &\leq& \frac{2}{p}\left(\frac{\theta+\tau}{2} B_{\eps}^0(\phi_{\eps}f^{p/2},\phi_{\eps}f^{p/2})+ (\eta (\theta)+\nu(\tau))\Vert f\Vert_{\eps,p}^p\right) \\
&\leq& \frac{(\theta + \tau)}{2}B_{\eps} (e^{a h}\phi_{\eps}f,e^{-a h}\phi_{\eps}f^{p-1}) \\
& & + \left( \frac{(\theta + \tau)a^2p}{4} +\frac{2}{p}(\eta(\theta) + \nu (\tau
))\right)  \Vert f\Vert_{\eps,p}^p.
\end{eqnarray*}
Letting now $\tau = \theta = \zeta/2>0$ yields
\begin{eqnarray*}
& &\EE_{\eps,p}(f) \\
&\leq& \frac{\zeta}{2}B_{\eps} (e^{a \widehat{h}}\phi_{\eps}f,e^{-a \widehat{h}}\phi_{\eps}f^{p-1})  + \left( \frac{\zeta a^2p }{4} +\frac{2}{p}(\eta(\zeta/2) + \nu (\zeta/2))\right)  \Vert f\Vert_{\eps,p}^p.
\end{eqnarray*}
Now, changing the constants suitably to $A,B>0$,
\begin{eqnarray*}
& & \eta(\zeta/2) + \nu (\zeta/2) \\
&=& k_1 + \frac{k_2\zeta}{2} -\frac{m+1}{4}\log(\zeta/2) + c_1 + \frac{c_2\zeta}{2} -\frac{1}{2}\log(\zeta/2)\\
&=& A +  B\zeta - \frac{m+3}{4}\log(\zeta),
\end{eqnarray*}
and therefore,
$$
\gamma (\zeta,p) = \frac{\zeta a^2p }{4} +\frac{2}{p}\left( A +  B\zeta -\frac{m+3}{4}\log(\zeta)\right).
$$
\hfill $\Box$

\subsection{The proof of the sub-Gaussian estimate}

In this final section, we will complete the proof of tightness by proving the uniform heat kernel bound stated in Proposition \ref{preparation}. \\

\noindent The operator $\widetilde{\HH_{\eps}}=\phi_{\eps}^{-1}(H_{\eps}-\lambda_{\eps})\phi_{\eps}$ is self - adjoint on the Hilbert space $L^2(L(1),\nu_{\eps})$, because 
$$
\widetilde{\HH_{\eps}} = u^{-1}\circ (H_{\eps}-\lambda_{\eps})\circ u,
$$
where $u: L^2(L(1),\nu_{\eps})\to L^2(L(1),\mu_{\Sa})$ is the unitary map given by $f\mapsto\phi_{\eps}f$. The domain of the operator is given by 
$$
\DD(\widetilde{\HH_{\eps}}) = u^{-1}\DD (H_{\eps}).
$$
By $\phi_{\eps}\geq 0$, the associated semigroup is still positivity-preserving and, by \cite{ReeSim:75}, Thm. X.55, p. 255, it determines contraction semigroups on $L^p(L(1),\phi_{\eps}^2\mu_{\Sa})$ for all $p\geq 1$.\\

\noindent Since the functions $\widehat{h}\in C^{\infty}(\overline{L(1)})$ are bounded and non-negative, multiplication with $e^{\pm a \widehat{h}}$ yields bounded and positivity-preserving endomorphisms for all spaces $L^p(L(1),\nu_{\eps})$ with $p\geq 1$. Thus, the operators 
$$
\HH_{\eps} = \frac{1}{2}\frac{e^{a \widehat{h}}}{\phi_{\eps}}(H_{\eps}-\lambda_{\eps})\frac{\phi_{\eps}}{e^{a \widehat{h}}},
$$
$\eps > 0$, are densely defined on $L^2(L(1),\nu_{\eps})$ and determine strongly continuous, positivity-preserving semigroups on all spaces $L^p(L(1),\nu_{\eps})$, $p\geq 1$ by 
$$
f\mapsto e^{a \widehat{h}} u^{-1}\left(e^{-\frac{t}{2}(H_{\eps}-\lambda_{\eps})}u(e^{-a \widehat{h}}f)\right) .
$$
Multiplication with $e^{\pm a\widehat{h}}$ leaves the domain $\DD (H_{\eps})= \bsob^1\cap\sob^2(L(1),\mu_{\Sa})$ invariant, hence the domain of the generator is
$$
\DD(\HH_{\eps}) = e^{-a\widehat{h}}\DD(\widetilde{\HH_{\eps}}) = u^{-1}(e^{-a\widehat{h}}\DD (H_{\eps})) = u^{-1}\DD (H_{\eps}) =  \DD(\widetilde{\HH_{\eps}}).
$$
From now on, we use the shorthand $f_s := e^{-s\HH_{\eps}}f$. By $\DD (\HH_{\eps})= u^{-1}(\bsob^1\cap\sob^2(L(1),\mu_{\Sa}))$, Theorem \ref{LP_logsobo} yields an estimate of the differential 
$$
\frac{d}{ds}\Vert f_s\Vert_{p(s)}^{p(s)} = \frac{B_{\eps}^0(\phi_{\eps}e^{a\widehat{h}}f_s, \phi_{\eps}e^{-a\widehat{h}}f_s^{p(s)-1}) + c(s) \EE_{\eps,p(s)}(f_s)}{\Vert f_s\Vert_{\eps,p(s)}^{p(s)}} \,\Vert f_s\Vert_{\eps,p(s)}
$$
for all $s>0$ whenever $f\in L^2(L(1),\nu_{\eps})$, $f\geq 0$. Here, $c(s) = \frac{d}{ds} \log p(s)$. The Gronwall - Lemma, together with a suitable choice of the function $p(s)$ yields thus an estimate of the form
$$
\Vert e^{-s\HH_{\eps}}u\Vert_{\infty}\leq e^M \Vert u\Vert_2
$$
for $f\in L^2(L(1),\nu_{\eps})$, $f\geq 0$. For general $f\in L^2(L(1),\nu_{\eps})$, the same estimate follows by
\begin{eqnarray*}
\Vert e^{-s\HH_{\eps}}f\Vert_{\infty} &=& \Vert \,\vert e^{-s\HH_{\eps}}f\vert\,\Vert_{\infty} = \Vert \,( e^{-s\HH_{\eps}}f )_+ +  (e^{-s\HH_{\eps}}f )_-\Vert_{\infty}\\
&=& \Vert \,e^{-s\HH_{\eps}}(u_+ + u_-)\Vert_{\infty} = \Vert \,e^{-s\HH_{\eps}}\vert u\vert\,\Vert_{\infty}\\
&\leq& e^M \Vert \,\vert f\vert\,\Vert_2 = e^M \Vert f\Vert_2.
\end{eqnarray*}

\noindent We are not going into the details of this proof. Instead, we use, with the same notations as above, the following statement, which is a direct consequence of \cite{Dav:89}, Theorem 2.2.7, p. 69. 

\begin{Theorem}\label{Davies} Let $\eps > 0$ be a continuous function defined for $2<p<\infty$ by
$$
\Gamma (p) = \gamma (\eps (p),p)
$$
If $t = \int_2^{\infty} p^{-1}\eps (p) dp$ and $M=\int_2^{\infty} p^{-1}\Gamma (p) dp$ are both finite, then the $C_0$-semigroup with generator 
$$
\HH_{\eps} = \frac{1}{2}\frac{e^{a \widehat{h}}}{\phi_{\eps}}(H_{\eps}-\lambda_{\eps})\frac{\phi_{\eps}}{e^{a\widehat{h}}},
$$ 
where $\widehat{h}\in C^{\infty}(\overline{L(1)})$ with $\sup_{W\in L(1)}\Vert d\widehat{h}\Vert_{\eps}\leq 1$, maps $L^2(L(1),\phi_{\eps}^2\mu_{\Sa})$ to $L^{\infty}(L(1),\mu_{\Sa})$ and
$$
\Vert   e^{-t\HH_{\eps}} \Vert_{\infty,2}\leq e^M.
$$
\end{Theorem}

\noindent We now come to the proof of Proposition \ref{preparation}. All that remains is to extract the estimate for the semigroup norm from the coefficients of the logarithmic Sobolev inequalities in Theorem \ref{LP_logsobo}. \\

\noindent{\bf Proof of Proposition \ref{preparation}.} Let $\Lambda > 1$ and $\zeta (p) := \Lambda 2^{-\Lambda}t p^{-\Lambda}$. Then
$$
\int_2^{\infty} p^{-1} \zeta (p) dp = t,
$$
and, for $\gamma (\zeta,p)$ from Theorem \ref{LP_logsobo}, we have
$$
p^{-1}\Gamma (p)  = \frac{\zeta(p)a^2 }{4} +\frac{2}{p^2}\left( A +  B\zeta(p) -\frac{m+3}{4}\log(\zeta(p)/t)-\frac{m+3}{4}\log(t)\right),
$$
and therefore
\begin{eqnarray*}
& & \int_2^{\infty} p^{-1} \Gamma (p) dp \\
&=& \frac{a^2 \Lambda t}{2(\Lambda - 1)} + A + \frac{Bt\Lambda}{\Lambda + 1} - \frac{m+3}{4} \log(t) + \frac{m+3}{4}(\log(\Lambda)-\Lambda)\\
&=& \frac{\Lambda}{2(\Lambda - 1)}a^2 t + C + Dt - \frac{m+3}{4} \log(t).
\end{eqnarray*}
That implies by Theorem \ref{Davies} for $t\leq 1$
\begin{equation}\label{mainest}
\Vert   e^{-t\HH_{\eps}} \Vert_{\infty,2}\leq e^{\frac{\Lambda}{2(\Lambda - 1)}a^2 t + C + Dt - \frac{m+3}{4} \log(t)} \leq Kt^{- \frac{m+3}{4}}e^{ \frac{\Lambda}{2(\Lambda - 1)}\rho^2 t}.
\end{equation}
As explained above, $\mathcal{H}_{\eps}$ is self-adjoint as an unbounded operator on the Hilbert space $L^2(L(1),\phi_{\eps}^2\mu_{\Sa})$. Hence, the adjoint of 
$$
e^{-t\HH_{\eps}}:L^2(L(1),\phi_{\eps}^2\mu_{\Sa})\to L^{\infty}(L(1),\phi_{\eps}^2\mu_{\Sa})
$$ 
is given by 
$$
\left(e^{-t\HH_{\eps}}\right)^* = e^{-t\mathcal{H}_{\eps}}:L^1(L(1),\phi_{\eps}^2\mu_{\Sa})\to L^2(L(1),\phi_{\eps}^2\mu_{\Sa})
$$ 
with $\Vert e^{-t\mathcal{H}_{\eps}}\Vert_{1,2}\leq\Vert e^{-t\mathcal{H}_{\eps}}\Vert_{2,\infty}$. That implies
$$
\Vert e^{-2t\mathcal{H}_{\eps}}\Vert_{\eps,1,\infty}= \Vert e^{-t\mathcal{H}_{\eps}}\circ\left(e^{-t\mathcal{H}_{\eps}}\right)^*\Vert_{\eps,1,\infty}\leq\Vert e^{-t\mathcal{H}_{\eps}}\Vert_{\eps,2,\infty}^2.
$$
\noindent By $\overline{h} = h\circ\pi$, $\pi(W)=x$, $\pi(W')=y$ and
$$
e^{-\frac{t}{2} H_{\eps}^0}f (W) = \int_{L(1)}K_t^{\eps}(W,W') f(W') \mu_{\Sa}(dW'),
$$
we obtain
\begin{eqnarray*}
& & e^{-t\HH_{\eps}}f (W) \\
&=& \int_{L(1)}e^{a (h(x)-h(y))}\frac{e^{\frac{t}{2}\left(\frac{\lambda_0}{\eps^2}-\lambda_{\eps}\right)}K_t^{\eps}(W,W')}{\phi_{\eps}(W)\phi_{\eps}(W')} f(y) \phi_{\eps}^2 (W')\mu_{\Sa}(dW').
\end{eqnarray*}
Therefore, the estimate (\ref{mainest}) implies the following estimate on the (symmetric and non-negative) kernel 
\begin{equation}\label{kernelesti}
K_t^{\eps}(W,W')\leq C_0 t^{- \frac{m+3}{2}}e^{ \frac{\Lambda}{\Lambda - 1}a^2 t}\phi_{\eps}(W)\phi_{\eps}(W')e^{\frac{t}{2}\left(\lambda_{\eps}-\frac{\lambda_0}{\eps^2}\right)}e^{a (h(y) - h(x))}.
\end{equation}
The sub-Gaussian estimate now follows from letting
$$
a :=\frac{h(x)-h(y)}{2kt},
$$
where $k:= \frac{\Lambda}{\Lambda -1} > 1$, since this implies
\begin{equation*}
K_t^{\eps}(W,W')\leq C_0 t^{- \frac{m+3}{2}}\phi_{\eps}(W)\phi_{\eps}(W')e^{\frac{t}{2}\left(\frac{\lambda_0}{\eps^2}-\lambda_{\eps}\right)}e^{-\frac{(h(y) - h(x))^2}{4kt}}.
\end{equation*}
Absorbing finally $e^{\frac{t}{2}\left(\frac{\lambda_0}{\eps^2}-\lambda_{\eps}\right)}$ into the constant by Proposition \ref{asymptotics}, we obtain some $C>0$ such that
\begin{equation}\label{subgaussian}
K_t^{\eps}(W,W')\leq C t^{- \frac{m+3}{2}}\phi_{\eps}(W)\phi_{\eps}(W') e^{-\frac{(h(y) - h(x))^2}{4kt}}.
\end{equation}
\hfill $\Box$

\end{document}